\documentclass[a4paper,12pt]{article}
\usepackage[utf8]{inputenc}
\usepackage{amsmath,amssymb,amsthm,mathtools}
\usepackage[margin=1.15in]{geometry}
\usepackage{tikz}
\usepackage[colorlinks=true,linkcolor=blue,citecolor=blue]{hyperref}
\usepackage{todonotes}

\usepackage{listings}
\newtheorem{thm}{Theorem}[section]
\newtheorem{prop}[thm]{Proposition}
\newtheorem{lem}[thm]{Lemma}
\newtheorem{cor}[thm]{Corollary}

\theoremstyle{definition}
\newtheorem{defn}[thm]{Definition}
\theoremstyle{remark}
\newtheorem{rem}[thm]{Remark}
\newtheorem{ex}[thm]{Example}                 

\DeclareMathOperator{\des}{des}
\DeclareMathOperator{\asc}{asc}
\DeclareMathOperator{\Des}{Des}
\DeclareMathOperator{\Ess}{Ess}
\DeclareMathOperator{\rk}{rk}
\DeclareMathOperator{\weak}{weak}
\DeclareMathOperator{\cen}{cen}
\DeclareMathOperator{\wt}{wt}
\DeclareMathOperator{\Wt}{Wt}
\DeclareMathOperator{\pk}{pk}
\newcommand{\N}{\mathbb{N}}
\newcommand{\Z}{\mathbb{Z}}
\newcommand{\Q}{\mathbb{Q}}
\newcommand{\cE}{\mathcal{E}}

\newcommand{\cL}{\mathcal{L}}
\newcommand{\cM}{\mathcal{M}}
\newcommand{\cF}{\mathcal{F}}

\newcommand{\one}{\mathbf{1}}
\newcommand{\Cat}{\mathrm{Cat}}
\newcommand{\id}{\mathrm{id}}

\title{Joint statistics on vexillary involutions and $\des$-Wilf classification of involutions avoiding a pattern of length four}
\date{}
\author{Niccol\`o Castronuovo\\ 
Universit\`a di Bologna, Campus di Cesena, 47521, ITALY \\
\texttt{niccolo.castronuovo2@unibo.it}\\
https://orcid.org/0000-0002-1054-0161}

\begin{document}
\maketitle

\begin{abstract}
We jointly enumerate various statistics over vexillary involutions among which there are the number of descents and the number of peaks. As a by-product, we get an equidistribution result between vexillary involutions and 3412-avoiding involutions, which, notably, implies that, for every $n$ and $k$, the number of involutions of length $n$
avoiding $3412$ and having $k$ descents equals the number of involutions of
length $n$ avoiding $2143$ and having $k$ ascents. The proof describes $2143$-avoiding involutions through Fulton's ranked
essential set, translates the count into a weighted enumeration of walks in the
quarter plane, and solves the latter explicitly.
These results, together with the determination of the generating functions of descents for involutions avoiding 3214 and 1243, respectively, obtained through the use of growth diagrams, allow us to complete the Wilf classification of the statistic $\des$ on every set of involutions avoiding a pattern of length four. Finally, we show that arbitrary
vexillary permutations are encoded by excursions of the same walk model used for involutions. 
\end{abstract}



\section{Introduction}

\subsection{Background}

The study of pattern-avoiding permutations has developed into a central area of enumerative combinatorics over the last three decades.  
Given a permutation $\pi\in S_n$ and a pattern $\tau\in S_k$, we say that $\pi$ \emph{avoids} $\tau$ if no subsequence of $\pi$ is order-isomorphic to $\tau$.  
The systematic study of such classes began with the foundational work of Knuth~\cite{Kn2}, and later expanded through the seminal contributions of Simion and Schmidt~\cite{Si}, among many others.  
Comprehensive accounts can be found, for instance, in the books by B\'ona~\cite{Bo} and Kitaev~\cite{Ki}.  

In parallel with the general theory of pattern-avoiding permutations, several authors have investigated \emph{involutions} avoiding specific patterns.  
Classical contributions include B\'ona et al.~\cite{Bonainv}, Guibert et al.~\cite{Gu,GuibertPergolaPinzani2001},  Barnabei et al.~\cite{Ba,Ba6}, and, more recently, Bean et al.~\cite{BeanGuttmannPantone2026}, among others.  
The avoidance of patterns by involutions often exhibits connections with lattice paths, continued fractions, and orthogonal polynomials, making enumeration strongly dependent on bijective or analytic techniques.

\emph{Vexillary permutations}, that is, permutations avoiding $2143$, were introduced
by Lascoux and Sch\"utzenberger in their study of Schubert polynomials, which
represent the classes of Schubert varieties in the cohomology ring of the complete
flag manifold \cite{LascouxSchutzenberger1982}. Vexillary permutations are precisely
those whose Schubert polynomial is a single flagged Schur function
\cite{LascouxSchutzenberger1982,Wachs1985}. They admit several equivalent
descriptions, for instance in terms of their diagrams and of their essential sets
\cite{Macdonald1991,Fulton1992,ErikssonLinusson1996}; Fulton's description by the
ranked essential set is the starting point of the present paper. From the enumerative
point of view, vexillary permutations are equinumerous with $1234$-avoiding
permutations \cite{BackelinWestXin2007}.

\emph{Vexillary involutions} were enumerated by Guibert, Pergola and Pinzani
\cite{GuibertPergolaPinzani2001} (see also \cite{Gu}), who proved that they are
counted by the Motzkin numbers; hence they are equinumerous with $3412$-avoiding
involutions, which are in bijection with Motzkin paths. 

When studying pattern avoidance, it is often the case that the enumeration is refined taking into accounts a statistic (see \cite{Ki}), e.g., the number of descents, the number of peaks, the major index, the number of right-to-left minima.

\subsection{Our contribution}

We refine this last result in the following sense. We give an explicit generating function for
vexillary involutions with respect to the number of descents, the number of peaks,
the numbers of initial and final fixed points and the lengths of the first and last
ascending runs (Theorem~\ref{thm:general}). As a consequence, we prove that the
descents of $2143$-avoiding involutions and the ascents of $3412$-avoiding involutions
are equidistributed, together with two further statistics (Theorem~\ref{thm:refined}).
Our approach encodes a vexillary involution, through its ranked essential set, by a
labeled walk in the quarter plane whose points record the numbers of dots north-west
and south-east of the essential boxes; all the statistics above can be read off this
walk, and the resulting weighted enumeration of walks turns out to have a solution in
product form. 
This result is the core of the $\des$-, $\asc$- and $\des$-$\asc$-Wilf
classifications of the patterns of length four on involutions, namely
the classification of the correspondences, for patterns
$\tau,\rho$ and every $n$ and $k$, between the number of
involutions in $I_n(\tau)$ with $k$ descents and the number of
involutions in $I_n(\rho)$ with $k$ descents (respectively ascents, or
the number of involutions in $I_n(\tau)$ with $k$ descents and the
number of involutions in $I_n(\rho)$ with $k$ ascents).
To get this classification we determine also the generating functions of descents (ascents, resp.) for involutions avoiding 3214 (1243, resp.), obtained through the use of growth diagrams. Finally, in an appendix, we show that arbitrary vexillary permutations are encoded by
excursions of the same walk model used for the involutions, and we derive a functional equation for the
corresponding generating function, whose explicit solution remains open.

\subsection{Outline} Section \ref{prel} introduces some basic notations and preliminary results and contains the statements of the main results. 
Section~\ref{sec:3412} deals with $I_n(3412)$, which is elementary, and proves the
$3412$ half of Theorem~\ref{thm:refined}. In Section~\ref{sec:ess} we use Fulton's
description of vexillary (that is, $2143$-avoiding) permutations by their ranked
essential sets \cite{Fulton1992} to encode $I_n(2143)$ by labeled  walks in $\N^2$
(Proposition~\ref{prop:walks}): the points of the walk record the numbers of dots
north-west and south-east of the essential boxes, and the statistics $\des$, $\pk$, $a$, $r$, $b$, $\bar r$ are read off the walk. In Section~\ref{sec:walks} we sum over the labels,
which expresses the generating function of Theorem~\ref{thm:general} in terms of a series
$W(r,s)$ of weighted walks; $W$ is computed in closed form in
Section~\ref{sec:solve}.
In Section~\ref{sec:wilf-blocks} we introduce some machinery (blocked involutions, their symmetric
fillings of Ferrers boards and the growth diagram bijections) necessary to deal with the 3214 and 1243 cases. We determine the descent-generating function for these cases in Sections~\ref{sec:wilf-3214} and~\ref{sec:wilf-1243}. In
Appendix~\ref{sec:vex} we show that arbitrary vexillary permutations correspond to
excursions of the same walk model, and we derive a functional equation for the corresponding generating function. Throughout the paper, all the identities stated to be checked with a computer algebra system are identities between formal power series; each of them is verified by the Mathematica code in Appendix~\ref{app:code}. 

\section{Preliminaries and statement of the main results}\label{prel}

\subsection{Preliminaries and notation}
Throughout, $\N=\{0,1,2,\dots\}$ and $[n]=\{1,\dots,n\}$. We denote by $S_n$
the set of permutations of $[n]$, written in one-line notation
$\sigma=\sigma(1)\sigma(2)\cdots\sigma(n)$, and by $I_n\subseteq S_n$ the set of
involutions; $S_0=I_0$ contains only the empty permutation. We write $\id_n$ for
the identity of $S_n$ and $\Omega_n=n(n-1)\cdots1.$ A
permutation $\sigma\in S_n$ \emph{contains} the pattern
$\tau\in S_k$ if there are positions $i_1<\dots<i_k$ such that
$\sigma(i_a)<\sigma(i_b)$ if and only if $\tau(a)<\tau(b)$, for all
$a,b\in[k]$; otherwise $\sigma$ \emph{avoids} $\tau$. We write $S_n(\tau)$ and
$I_n(\tau)$ for the sets of permutations and of involutions of $[n]$ avoiding
$\tau$. Set $I(\tau):=\cup_{n\geq 0}I_n(\tau).$
For $\sigma\in S_n$,
\begin{gather*}
\Des(\sigma)=\{i\in[n-1]:\sigma(i)>\sigma(i+1)\},\qquad
\operatorname{Asc}(\sigma)=\{i\in[n-1]:\sigma(i)<\sigma(i+1)\},\\
\des(\sigma)=|\Des(\sigma)|,\qquad \asc(\sigma)=|\operatorname{Asc}(\sigma)|,
\end{gather*}
so that $\asc(\sigma)=n-1-\des(\sigma)$ if $n\ge1$, and $\des=\asc=0$ on $S_0$.
For $\sigma\in S_n$ let $$a(\sigma)=\max\{a\in\{0,\dots,n\}:\sigma(i)=i\text{ for all }i\le a\}$$
be the number of \textit{initial fixed points} of $\sigma$ and, if $\sigma\neq\id_n,$ denote by $f(\sigma)$ the \textit{length of the first increasing run}, namely
let $$f(\sigma)=\min\Des(\sigma)\qquad\text{and set}\qquad r(\sigma)=f(\sigma)-a(\sigma).$$
Let also $$b(\sigma)=\max\{b\in\{0,\dots,n\}:\sigma(i)=i\text{ for all }i>n-b\}$$ be the
number of \textit{final fixed points} of $\sigma$ and, if $\sigma\neq\id_n$, let
$$\ell(\sigma)=n-\max\Des(\sigma)$$ be the \textit{length of the last increasing run} of $\sigma$
and set $$\bar r(\sigma)=\ell(\sigma)-b(\sigma).$$ Finally, let
$\pk(\sigma)=|\{i\in\{2,\ldots,n-1\}:\sigma(i-1)<\sigma(i)>\sigma(i+1)\}|$ be the \textit{number of
peaks} of $\sigma$.

\begin{ex}\label{ex:stats}
Let $\sigma=1\,5\,7\,4\,2\,6\,3\,8$, an involution avoiding $2143$. Its ascending runs
are $157$, $4$, $26$, $38$, and $\Des(\sigma)=\{3,4,6\}$, so $\des(\sigma)=3$. We have
$a(\sigma)=1$, $f(\sigma)=3$, $r(\sigma)=2$, $b(\sigma)=1$, $\ell(\sigma)=8-6=2$ and
$\bar r(\sigma)=1$. The peaks of $\sigma$ are $3$ and $6$, so $\pk(\sigma)=2$; the
descent $4$ is not a peak, because the run ending at $4$ has length $1$.
\end{ex}

A \emph{Motzkin path} of length $n$ is a word $w=w_1\cdots w_n$ over
$\{U,H,D\}$ such that every prefix contains at least as many $U$'s as $D$'s and
$w$ contains as many $U$'s as $D$'s. Let $\cM_n$ be the set of Motzkin paths of
length $n$, $\cM=\bigcup_{n\ge0}\cM_n$, and let $|w|$ denote the length of $w$.
A \emph{center} of $w\in\cM$ is either a letter $H$ of $w$ or a factor $UD$ of $w$;
for $w\neq\emptyset$ we denote by $\cen(w)$ the number of centers of $w,$ and we set $\cen(\emptyset)=0$. The \emph{tower} of a center is the longest
factor of $w$ of the form $U^jHD^j$ (if the center is an $H$) or $U^jD^j$ with
$j\ge1$ (if the center is a factor $UD$) having that center in the middle. A \emph{weak valley} of a Motzkin path $w$ is a factor of $w$ belonging to
$\{HU,\,HH,\,DU,\,DH\}$ and $\weak(w)$ is the number of weak valleys of $w.$ 

\begin{rem}\label{weak_cen}
For every nonempty
$w\in\cM$, $\weak(w)=\cen(w)-1.$ 
\begin{proof}
Cut $w$ at its weak valleys: $w=w^{(1)}\cdots w^{(\ell)}$ with
$\ell=\weak(w)+1$, each $w^{(i)}$ containing no weak valley. In $w^{(i)}$ a letter
$H$ or $D$ can only be followed by $D$, so $w^{(i)}=U^{e}H^{\epsilon}D^{g}$ with
$\epsilon\in\{0,1\}$. Every $w^{(i)}$ begins with $U$ or $H$ (because $w$ does,
and a weak valley ends with $U$ or $H$) and ends with $H$ or $D$ (because $w$
does, and a weak valley begins with $H$ or $D$). Hence each $w^{(i)}$ contains
exactly one center: its letter $H$ if $\epsilon=1$, its factor $UD$ if
$\epsilon=0$. No center meets two consecutive factors, since $UD$ is not a weak
valley. Therefore $w$ has exactly $\ell=\weak(w)+1$ centers.
\end{proof}
\end{rem}

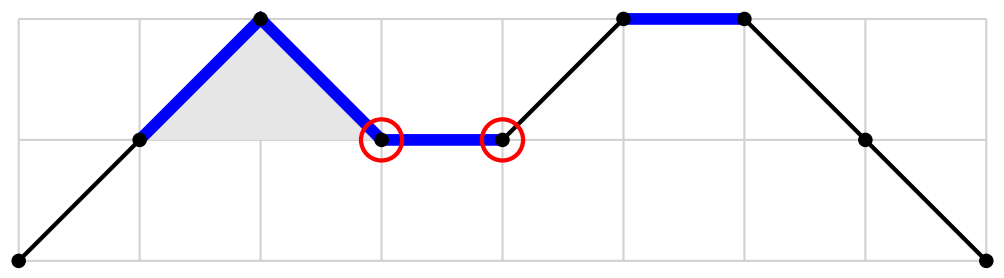
\begin{figure}[ht]
\centering
\begin{tikzpicture}[scale=0.8]
  \draw[gray!35] (0,0) grid (8,2);
  \fill[gray!20] (1,1) -- (2,2) -- (3,1) -- cycle;
  \draw[thick] (0,0) -- (1,1) -- (2,2) -- (3,1) -- (4,1) -- (5,2) -- (6,2) -- (7,1) -- (8,0);
  \draw[line width=2.2pt,blue] (1,1) -- (2,2) -- (3,1);
  \draw[line width=2.2pt,blue] (3,1) -- (4,1);
  \draw[line width=2.2pt,blue] (5,2) -- (6,2);
  \foreach \p in {(0,0),(1,1),(2,2),(3,1),(4,1),(5,2),(6,2),(7,1),(8,0)} \fill \p circle (0.06);
  \draw[red,thick] (3,1) circle (0.17);
  \draw[red,thick] (4,1) circle (0.17);
  \foreach \i/\l in {1/U,2/U,3/D,4/H,5/U,6/H,7/D,8/D} \node[below] at ({\i-0.5},-0.15) {\scriptsize $\l$};
\end{tikzpicture}
\caption{The Motzkin path $w=UUDHUHDD$ of Example~\ref{ex:motzkin}. As usual, we
represent a Motzkin path in the plane, the letters $U$, $H$ and $D$ corresponding to
the steps $(1,1)$, $(1,0)$ and $(1,-1)$. The centers are drawn in bold; the circles
mark the middle points of the two weak valleys $DH$ and $HU$; the shaded triangle is
the tower of the leftmost center.}
\label{fig:motzkin}
\end{figure}

\begin{ex}\label{ex:motzkin}
Let $w=UUDHUHDD\in\cM_8$ (Figure~\ref{fig:motzkin}). Its centers are the factor $UD$
formed by its second and third letters and the letters $H$ in positions $4$ and $6$,
so $\cen(w)=3$. Its weak valleys are the factors $DH$ (positions $3,4$) and $HU$
(positions $4,5$), so $\weak(w)=2=\cen(w)-1$, in accordance with
Remark~\ref{weak_cen}. The towers of the three centers are $UD$ (positions $2,3$),
$H$ (position $4$) and $UHD$ (positions $5$--$7$).
\end{ex}

For $\pi\in I_n$ let $\beta(\pi)=w_1\cdots w_n$, where $w_i=U$, $H$ or $D$
according as $\pi(i)>i$, $\pi(i)=i$ or $\pi(i)<i$.

\begin{lem}\label{lem:asc3412}
The map $\beta$ restricts to a bijection $\beta:I(3412)\to \mathcal M$ such that
$\beta(I_n(3412))=\mathcal M_n$ for every $n\ge0$, and
$\asc(\pi)=\cen(\beta(\pi))-1$ for every $\pi\in I_n(3412)$ with $n>0$.
\end{lem}

\begin{proof}
By \cite[Lemma 3.2]{dmtcs:5175}, $\beta$ maps $I_n(3412)$ bijectively onto $\cM_n$,
and $i\in[n-1]$ is an ascent of $\pi$ if and only if $w_iw_{i+1}$ is a weak valley of $w=\beta(\pi).$ Hence $\asc(\pi)$ is the number $\weak(w)$ of weak valleys
of $w.$ This concludes the proof since $\weak(w)=\cen(w)-1,$ as observed in Remark \ref{weak_cen}. 
\end{proof}

If $\pi\in I_n(3412)$ with
$n\ge1$, then $\beta(\pi)$ is a nonempty Motzkin path (Lemma~\ref{lem:asc3412});
we let $T(\pi)$ be the tower of its leftmost center, and we denote by $s(\pi)$
the length of $T(\pi)$ and by $u(\pi)$ the number of letters $U$ of $\beta(\pi)$
to the left of $T(\pi)$.

\begin{ex}
For $\pi=6\,3\,2\,5\,4\,1$ we have $\beta(\pi)=UUDUDD$, with centers $UD$ in
positions $2,3$ and $4,5$, so that $\cen(\beta(\pi))=2$ and
$\asc(\pi)=\cen(\beta(\pi))-1=1$. The first tower is
the factor $UD$ in positions $2,3$, so $s(\pi)=2$ and $u(\pi)=1$.
\end{ex}

All generating functions are elements of the ring $R=\Q[[x,t,\omega,v,w,\bar v,\bar w]]$.
We fix once and for all the elements
\begin{equation}\label{eq:zyc}
z=\frac{x}{1-x},\qquad y=\frac{x^2}{1-x^2},\qquad c=tz=\frac{tx}{1-x}
\end{equation}
of $R.$ Put $q_\omega=1-t(1-\omega)x$ and
$B_\omega=1+x^2-tx-t\omega x^2-t(1-\omega)x^3$, and let $\theta_\omega$ be the unique
element of $R$ with zero constant term such that
\begin{equation}\label{eq:starom}
x\,q_\omega\,\theta_\omega^2-B_\omega\,\theta_\omega+x\,q_\omega=0,
\qquad\text{namely}\qquad
\theta_\omega=\frac{B_\omega-\sqrt{B_\omega^2-4x^2q_\omega^2}}{2x\,q_\omega}.
\end{equation}

When $\omega=1$ we will not indicate the subscript. So, for example, $\theta=\theta_1$
is the unique element of $R$ with zero constant term such
that
\begin{equation}\label{eq:star}
x\,\theta^2-\bigl(1+x^2-tx(1+x)\bigr)\,\theta+x=0 ,
\end{equation}
namely
\[
\theta=\frac{1+x^2-tx(1+x)-\sqrt{\bigl(1+x^2-tx(1+x)\bigr)^2-4x^2}}{2x}.
\]

For a proposition $\mathcal{P}$, set $$[\mathcal{P}]=\begin{cases}
    1&\mbox{if }\mathcal{P}\mbox{ holds}\\
    0&\mbox{otherwise}
\end{cases}$$ 
(the previous notation is not confused with the squared brackets representing sets, as in $[n],$ since in this case the content of the brackets is an integer, not a proposition).

\subsection{Results}
Now we state our main result.

\begin{thm}\label{thm:general}
We have,
\begin{multline*}
\sum_{n\ge1}\ \sum_{\sigma\in I_n(2143)\setminus\{\id_n\}}
t^{\des\sigma}\omega^{\pk\sigma}v^{a(\sigma)}w^{r(\sigma)}\bar v^{\,b(\sigma)}\bar w^{\,\bar r(\sigma)}x^n
\\
=\frac{t\,x\,\theta_\omega\,w\bar w\,(1-\bar w x^2)}{q_\omega\,(1-x\bar v)(1-x\bar w\theta_\omega)}
\Bigl(1-\omega+\frac{\omega\,(1-wx^2)}{(1-xv)(1-xw\theta_\omega)(1-w\bar w x^2)}\Bigr).
\end{multline*}
\end{thm}

For $\omega=\bar v=\bar w=1$ the right-hand side of Theorem~\ref{thm:general} becomes
$\frac{t\,x\,(1+x)\,\theta\,w}{(1-x\theta)(1-xv)(1-xw\theta)}=\frac{w(\theta-x)}{(1-xv)(1-xw\theta)}$,
by \eqref{eq:star}. After these substitutions, we can compare such generating function with another taking into account statistics over $3412$-avoiding involutions.

\begin{thm}\label{thm:refined}
For every $n\ge1$, the triple $(\des,a,r)$ on $I_n(2143)\setminus\{\id_n\}$ has
the same distribution as the triple $(\asc,s-1,u+1)$ on
$I_n(3412)\setminus\{\Omega_n\}$. More precisely, both
generating functions
\[
\sum_{n\ge1}\ \sum_{\sigma\in I_n(2143)\setminus\{\id_n\}}t^{\des\sigma}v^{a(\sigma)}w^{r(\sigma)}x^n
\quad\text{and}\quad
\sum_{n\ge1}\ \sum_{\pi\in I_n(3412)\setminus\{\Omega_n\}}t^{\asc\pi}v^{s(\pi)-1}w^{u(\pi)+1}x^n
\]
are equal to
\[
\frac{w\,(\theta-x)}{(1-xv)(1-xw\theta)} .
\]
\end{thm}

\begin{cor}\label{cor:peaks}
Let $F_\omega=\sum_{n\ge0}\sum_{\sigma\in I_n(2143)}t^{\des\sigma}\omega^{\pk\sigma}x^n$. Then
$$F_\omega = 1 + (\theta_\omega-x)/(t \omega x).$$
\end{cor}

Set
\[
F(x,t)=\sum_{n\ge 0}\ \sum_{\sigma\in I_n(2143)}t^{\des(\sigma)}x^n,
\qquad
N(x,t)=\sum_{n\ge 0}\ \sum_{\pi\in I_n(3412)}t^{\asc(\pi)}x^n,\]
\[O(x,t)=\sum_{n\ge 0}\ \sum_{\pi\in I_n(3214)}t^{\des(\pi)}x^n,\qquad
P(x,t)=\sum_{n\ge 0}\ \sum_{\pi\in I_n(1243)}t^{\asc(\pi)}x^n.
\]

\begin{thm}\label{wilf_equidist_des}
We have
\begin{align*}
F(x,t)&=N(x,t)=O(x,t)=P(x,t)=\frac{1+\theta}{1-x\theta}\\
&=\frac{1}{2tx^2}\Bigl(1-x^2+x^2t-xt-\sqrt{x^4(t-1)^2+x^3(2t^2-2t)+x^2(t^2-2t-2)-2tx+1}\Bigr).
\end{align*}
As a consequence,
for all $n\ge 1$ and $k\ge 0$,
\[
\bigl|\{\sigma\in I_n(2143):\des(\sigma)=k\}\bigr|=\bigl|\{\sigma\in I_n(3412):\asc(\sigma)=k\}\bigr|
=\bigl|\{\sigma\in I_n(3214):\des(\sigma)=k\}\bigr|\]
\[=\bigl|\{\sigma\in I_n(1243):\asc(\sigma)=k\}\bigr|=[x^n]\; z^{k+1}\sum_{j=0}^{k}\binom{k+j}{k-j}\,\Cat(j)\,y^j,
\]
where $\Cat(j)=\frac{1}{j+1}\binom{2j}{j}$ and $[x^n]$ denotes the coefficient of $x^n$.
It is also the number of Motzkin paths of
length $n$ with $k$ weak valleys, that is, with $k+1$ centers.
\end{thm}

Setting $t=1$ in the expression of $F(x,t)$ one recovers the theorem of Guibert,
Pergola and Pinzani~\cite{GuibertPergolaPinzani2001} that $|I_n(2143)|$ is the
$n$-th Motzkin number.

\section{Involutions avoiding \texorpdfstring{$3412$}{3412}}\label{sec:3412}

\begin{prop}\label{prop:N}
The series $N$ satisfies
\begin{equation}\label{eq:motzkin}
N=1+x\bigl(1+t(N-1)\bigr)+x^2N\bigl(1+t(N-1)\bigr).
\end{equation}
Moreover, we have $\theta=x\bigl(1+t(N-1)\bigr)$, where $\theta$ is the series defined
by \eqref{eq:star}, and 
\begin{equation}\label{eq:Ncat}
y\widetilde N^2-(1-c)\widetilde N+c=0,
\end{equation}
where $\widetilde N=t(N-1).$
\end{prop}

\begin{proof}
By Lemma~\ref{lem:asc3412},
$N=1+\sum_{w\in\cM\setminus\{\emptyset\}}t^{\cen(w)-1}x^{|w|}$, that is,
$\sum_{w\in\cM}t^{\cen(w)}x^{|w|}=1+t(N-1)$. A nonempty Motzkin path begins with
$U$ or $H$ and ends with $H$ or $D$, and it can be written uniquely either as
$w=Hw'$ or as $w=Uw''Dw'$ with $w',w''\in\cM$. The empty path contributes $1$.

If $w=Hw'$, then $\cen(w)=1+\cen(w')$, since the centers of $w$ are its first
letter and the centers of $w'$ (no factor $UD$ contains the first letter $H$).
These words contribute $x\sum_{w'\in\cM}t^{\cen(w')}x^{|w'|}=x\bigl(1+t(N-1)\bigr)$
to $N$.

If $w=Uw''Dw'$, put $v=Uw''D$. If $w''$ is empty then $\cen(v)=1$, the only
center of $v=UD$ being $v$ itself. Otherwise the first two letters of $v$ are
$UU$ or $UH$ and the last two are $HD$ or $DD$, so the centers of $v$ are those
of $w''$ and $\cen(v)=\cen(w'')$. Hence
$\sum_{w''\in\cM}t^{\cen(v)-1}x^{|v|}=x^2N$. Finally $\cen(w)=\cen(v)+\cen(w')$,
since $v$ ends with $D$, so that no factor $UD$ meets both $v$ and $w'$. These
words contribute $x^2N\bigl(1+t(N-1)\bigr)$. This proves
\eqref{eq:motzkin}.

The remaining assertions are trivial algebraic manipulations, verified by a
symbolic algebra program. 
\end{proof}

\begin{prop}\label{prop:G3412}
We have
\[
\sum_{n\ge1}\ \sum_{\pi\in I_n(3412)\setminus\{\Omega_n\}}t^{\asc\pi}v^{s(\pi)-1}w^{u(\pi)+1}x^n
=\frac{w\,(\theta-x)}{(1-xv)(1-xw\theta)} .
\]
\end{prop}

\begin{proof}
Let $\pi\in I_n(3412)$, $n\ge1$, and $w=\beta(\pi)$. By Lemma~\ref{lem:asc3412},
$\asc(\pi)=\cen(w)-1$ is the number of centers of $w$ minus one.
The prefix of $w$ preceding its leftmost center contains neither $H$ nor $UD$,
and $w$ does not begin with $D$; hence either $w=U^mH\cdots$ with $m\ge0$ (center
$H$), or $w=U^mD\cdots$ with $m\ge1$ (center $UD$, formed by the $m$-th letter
and the following $D$). Let $k$ be the number of letters $D$ immediately
following the center ($H$) or its letter $U$ (center $UD$). Then $k\le m$,
$T(\pi)=U^kHD^k$ or $U^kD^k$, so that $s(\pi)=2k+1$ or $2k$, $u(\pi)=m-k$, and
\[
w=U^{m-k}\,T(\pi)\,X,
\]
where $X$ goes from height $h=m-k$ to height $0$ and does not begin with $D$.
Such an $X$ factors uniquely as $X=E_0DE_1D\cdots DE_h$ with $E_i\in\cM$ and
$E_0\neq\emptyset$ if $h\ge1$. The centers of $w$ are the center of $T(\pi)$ and
the centers of the $E_i$: a factor $UD$ cannot straddle two consecutive pieces,
because each piece preceding a $D$ of the factorization (namely $T(\pi)$ or some
$E_i$) is empty or ends with $H$ or $D$. Let $M=\sum_{w\in\cM}t^{\cen(w)}x^{|w|}$
be the generating function of Motzkin paths with weight $t$ per center; by
Lemma~\ref{lem:asc3412} and Proposition~\ref{prop:N}, $M=1+t(N-1)=\theta/x$. Then
the paths with given $k$ and $h$ contribute
\[
\bigl(x^{2k+1}v^{2k}\ \text{or}\ x^{2k}v^{2k-1}\bigr)\cdot([h=0]w\,M+[h\geq 1]w^{h+1}x^{2h}(M-1)M^h).
\]
The excluded involution $\Omega_n$ corresponds to $h=0$ and $X=\emptyset$. Summing
over $k$ gives the factor $\frac{x+x^2v}{1-x^2v^2}=\frac{x}{1-xv}$, and summing
over $h$ gives
\[
wM-w+\sum_{h\ge1}w^{h+1}x^{2h}(M-1)M^h=\frac{w(M-1)}{1-x^2wM}
=\frac{w(\theta-x)}{x(1-xw\theta)} .\qedhere
\]
\end{proof}

\section{The ranked essential set of a vexillary involution}\label{sec:ess}

\subsection{Diagrams, essential sets, ranks}
\begin{figure}[ht]
\centering
\begin{tikzpicture}[scale=0.62]
  \def\n{8}
  \fill[black!18] (0,0) rectangle (\n,\n);
  \foreach \i/\j in {1/1,1/2,1/3,1/4,1/5,1/6,2/1,3/1,3/3,3/4,3/5,4/1,4/3,
                     5/1,5/3,5/5,6/1}
    \fill[white] ({\j-1},{\n-\i}) rectangle ({\j},{\n-\i+1});
  \draw[black!45] (0,0) grid (\n,\n);
  \draw (0,0) rectangle (\n,\n);
  \foreach \i/\j in {1/7,2/2,3/6,4/4,5/8,6/3,7/1,8/5}
    \fill ({\j-0.5},{\n-\i+0.5}) circle (0.14);
  \foreach \i/\j/\r in {1/6/0,3/5/1,5/5/2,5/3/1,6/1/0} {
    \draw[line width=1.3pt] ({\j-1},{\n-\i}) rectangle ({\j},{\n-\i+1});
    \node at ({\j-0.5},{\n-\i+0.5}) {\small$\r$};
  }
  \foreach \i in {1,...,\n} {
    \node[left]  at (0,{\n-\i+0.5}) {\scriptsize$\i$};
    \node[above] at ({\i-0.5},\n)   {\scriptsize$\i$};
  }
\end{tikzpicture}
\caption{The involution $\sigma=7\,2\,6\,4\,8\,3\,1\,5$. The dots are the
boxes $(i,\sigma(i))$; the white boxes form the diagram $D(\sigma)$; the
boxes with a thick border form the essential set $\Ess(\sigma)$, and the
number inside each of them is the value of $\rk_\sigma$.}
\label{fig:essential}
\end{figure}
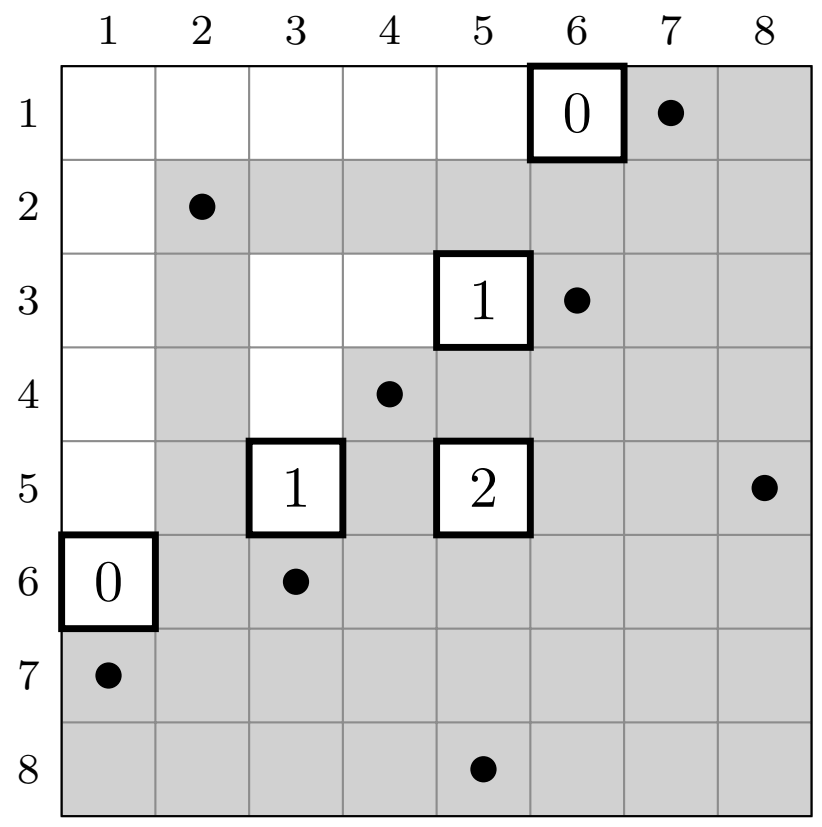

Let $\sigma\in S_n$. Following \cite{Fulton1992,Macdonald1991}, the
\emph{diagram}, the \emph{essential set} and the \emph{rank function} of
$\sigma$ are
\begin{align*}
D(\sigma)&=\{(i,j)\in[n]^2:\ \sigma(i)>j\ \text{and}\ \sigma^{-1}(j)>i\},\\
\Ess(\sigma)&=\{(i,j)\in[n-1]^2:\ \sigma(i)>j,\ \ \sigma(i+1)\le j,\ \
\sigma^{-1}(j)>i,\ \ \sigma^{-1}(j+1)\le i\},\\
\rk_\sigma(i,j)&=\bigl|\{h\in[i]:\ \sigma(h)\le j\}\bigr|\qquad(i,j\in[n]).
\end{align*}
We draw $[n]^2$ as an $n\times n$ array of boxes, with rows numbered from
top to bottom and columns from left to right, and we put a dot in the boxes
$(i,\sigma(i))$, $i\in[n]$. Then $D(\sigma)$ is the set of boxes that
remain after deleting, for every dot, the box containing it and all the
boxes to its right in the same row or below it in the same column, and
$\rk_\sigma(i,j)$ is the number of dots weakly north-west of the box
$(i,j)$; see Figure~\ref{fig:essential}. Note that
$D(\sigma)\subseteq[n-1]^2$.

The definition of $\Ess(\sigma)$ is \cite[(3.8)]{Fulton1992}. The essential
set is the set of south-east corners of $D(\sigma)$, that is,
\begin{equation}\label{eq:corners}
\Ess(\sigma)=\{(i,j)\in D(\sigma):\ (i+1,j)\notin D(\sigma)\ \text{and}\
(i,j+1)\notin D(\sigma)\},
\end{equation}
as in \cite[Section~2]{ErikssonLinusson1996}. Indeed, if
$(i,j)\in D(\sigma)$, then $(i+1,j)\notin D(\sigma)$ if and only if
$\sigma(i+1)\le j$ and $(i,j+1)\notin D(\sigma)$ if and only if
$\sigma^{-1}(j+1)\le i.$ To prove the first equivalence notice that, because \(\sigma\) is a permutation, we cannot have both $\sigma(i+1)=j$ and $\sigma^{-1}(j)>i+1,$
hence, if \(\sigma(i+1)>j\), then \(\sigma^{-1}(j)\neq i+1\). Together with \(\sigma^{-1}(j)>i\), this gives $\sigma^{-1}(j)>i+1.$ The proof of the second equivalence is symmetric.

Moreover, $\sigma(i+1)\le j$ implies
$(i+1,j+1)\notin D(\sigma)$; hence \eqref{eq:corners} also agrees with the
description of $\Ess(\sigma)$ given in \cite[p.~392]{Fulton1992}, where the
(redundant) condition $(i+1,j+1)\notin D(\sigma)$ is added.

The \emph{ranked essential set} of $\sigma$ is
\[
\cE(\sigma)=\bigl\{(p,q,\rk_\sigma(p,q)):\ (p,q)\in\Ess(\sigma)\bigr\}.
\]
A permutation is \emph{vexillary} if it avoids $2143$ (see, e.g., \cite{Macdonald1991,GuibertPergolaPinzani2001}).

\begin{thm}[Eriksson--Linusson \cite{ErikssonLinusson1996},Fulton \cite{Fulton1992}]\label{thm:fulton}
Let $n\ge 1$.
\begin{enumerate}
\item[(a)] A permutation $\sigma\in S_n$ is determined by $\cE(\sigma)$.
\item[(b)] Let $k\ge0$ and let $(p_i,q_i,\rho_i)$, $1\le i\le k$, be triples
of (nonnegative) integers such that
\begin{gather}
p_1\le p_2\le\dots\le p_k,\qquad q_1\ge q_2\ge\dots\ge q_k,\label{eq:F1}\\
0<p_1-\rho_1<p_2-\rho_2<\dots<p_k-\rho_k,\qquad
q_1-\rho_1>q_2-\rho_2>\dots>q_k-\rho_k>0,\label{eq:F2}\\
\rho_i\ge0\quad\text{and}\quad p_i+q_i-\rho_i\le n\qquad(1\le i\le k).\label{eq:F3}
\end{gather}
Then there is a unique $\sigma\in S_n$ such that
$\cE(\sigma)=\{(p_i,q_i,\rho_i):1\le i\le k\}$, and $\sigma$ is vexillary.
Conversely, for every vexillary $\sigma\in S_n$ the elements of $\cE(\sigma)$
can be indexed, in a unique way, so that \eqref{eq:F1}, \eqref{eq:F2} and
\eqref{eq:F3} hold.
\end{enumerate}
\end{thm}

\begin{proof}
Part (a) is \cite[Lemma~3.10(b)]{Fulton1992}. Part (b) is a reformulation of
\cite[Proposition~2.2]{ErikssonLinusson1996}, which states that a set
$E\subseteq[n]^2$ with a function $r\colon E\to\mathbb{Z}$ is the essential
set, with its rank function, of a vexillary permutation in $S_n$ if and only
if
\begin{enumerate}
\item[(EL1)] $\min\{i,j\}>r(i,j)\ge0$ and $r(i,j)+n\ge i+j$ for all $(i,j)\in E$;
\item[(EL2)] $i-i'>r(i,j)-r(i',j')>j-j'$ for all distinct $(i,j),(i',j')\in E$
such that $i\ge i'$, $j\le j'$ and
$E\cap\bigl([i',i]\times[j,j']\bigr)=\{(i,j),(i',j')\}$;
\item[(EL3)] there are no $(i_1,j_1),(i_2,j_2)\in E$ with $i_1<i_2$ and $j_1<j_2$.
\end{enumerate}
Condition (EL3) says that the elements of $E$ can be listed as
$(p_1,q_1),\dots,(p_k,q_k)$ so that \eqref{eq:F1} holds. With this listing
and $\rho_i=r(p_i,q_i)$, the pairs to which (EL2) applies are the pairs
$\{(p_i,q_i),(p_{i+1},q_{i+1})\}$, and (EL2) reads
$p_i-\rho_i<p_{i+1}-\rho_{i+1}$ and $q_i-\rho_i>q_{i+1}-\rho_{i+1}$; finally
(EL1) reads $p_i-\rho_i>0$, $q_i-\rho_i>0$, $\rho_i\ge0$ and
$p_i+q_i-\rho_i\le n$. Hence (EL1), (EL2) and (EL3) together are equivalent to
the existence of an indexing satisfying \eqref{eq:F1}--\eqref{eq:F3}; this
indexing is unique because the numbers $p_i-\rho_i$ are strictly increasing.
The uniqueness of $\sigma$ follows from (a).
\end{proof}

\begin{rem}\label{rem:fulton}
The conditions \eqref{eq:F1} and \eqref{eq:F2} are those of
\cite[Proposition~9.6]{Fulton1992}, where the essential boxes are written
$(b_i,a_i)$ (row, column) with ranks $r_i$ and indexed in the reverse order,
$(p_i,q_i,\rho_i)=(b_{k+1-i},a_{k+1-i},r_{k+1-i})$. Fulton proves the
existence of $\sigma$ in $S_n$ for $n\ge a_k+b_1$; the sharp condition
$p_i+q_i-\rho_i\le n$ in \eqref{eq:F3} is due to Eriksson and Linusson.
\end{rem}

\begin{lem}\label{lem:inv-des}
Let $\sigma\in S_n$.
\begin{enumerate}
\item[(a)] $\cE(\sigma^{-1})=\{(q,p,\rho):\ (p,q,\rho)\in\cE(\sigma)\}$.
\item[(b)] $\Des(\sigma)=\{p:\ (p,q)\in\Ess(\sigma)\text{ for some }q\}$.
\end{enumerate}
\end{lem}

\begin{proof}
(a) The diagram $D(\sigma^{-1})$ is the transpose of $D(\sigma)$
\cite[(1.21)(i)]{Macdonald1991}; hence, by \eqref{eq:corners},
$\Ess(\sigma^{-1})$ is the transpose of $\Ess(\sigma)$. Moreover
$\rk_{\sigma^{-1}}(q,p)=\rk_\sigma(p,q)$, since both numbers count the pairs
$(h,l)\in[p]\times[q]$ with $\sigma(h)=l$.

(b) See \cite[proof of Proposition~9.18]{Fulton1992} or
\cite[Section~5]{ErikssonLinusson1996}.
\end{proof}

\begin{ex}\label{ex:essential}
Let $n=8$ and $\sigma=7\,2\,6\,4\,8\,3\,1\,5\in I_8(2143)$. The diagram
$D(\sigma)$ consists of the $17$ white boxes of
Figure~\ref{fig:essential}, and
\[
\cE(\sigma)=\{(1,6,0),\ (3,5,1),\ (5,5,2),\ (5,3,1),\ (6,1,0)\}.
\]
The triples are listed in the order of Theorem~\ref{thm:fulton}(b): the
sequences $(p_i-\rho_i)=(1,2,3,4,6)$ and $(q_i-\rho_i)=(6,4,3,2,1)$ are
strictly monotone, and $(p_i+q_i-\rho_i)=(7,7,8,7,7)$. Since $\sigma$ is an
involution, $\cE(\sigma)$ is symmetric, in accordance with
Lemma~\ref{lem:inv-des}(a); the rows of the essential boxes are
$1,3,5,6$, and indeed $\Des(\sigma)=\{1,3,5,6\}$.
\end{ex}

\subsection{The encoding by labeled  walks}

For a box $(p,q)\in D(\sigma)$, the $n$ dots of $\sigma$ split into the dots
north-west of $(p,q)$ (rows $\le p$, columns $\le q$), north-east (rows $\le p$,
columns $>q$), south-west and south-east (rows $>p$, columns $>q$). Their
numbers are, respectively,
\[
\rk_\sigma(p,q),\qquad p-\rk_\sigma(p,q),\qquad q-\rk_\sigma(p,q),\qquad
n-p-q+\rk_\sigma(p,q).
\]
Conditions \eqref{eq:F2} and \eqref{eq:F3} say that, along the essential set, the
number of dots north-east increases strictly, the number of dots south-west
decreases strictly, and the numbers of dots north-west and south-east are
nonnegative. This suggests the following encoding.

\begin{defn}\label{def:walks}
A \emph{labeled  walk} is a finite sequence
$\lambda=\bigl((u_i,\alpha_i,\beta_i)\bigr)_{i=1}^m$, $m\ge0$, of triples of
integers such that, setting $(r_0,s_0)=(0,0)$ and
\[
r_i=r_{i-1}+\alpha_i-u_i,\qquad s_i=s_{i-1}+\beta_i-u_i\qquad(1\le i\le m),
\]
the following hold:
\begin{enumerate}
\item[(S1)] $u_i\ge1$, $\alpha_i\ge0$, $\beta_i\ge0$ and
$\alpha_i+\beta_i\ge u_i+1$ for $1\le i\le m$;
\item[(S2)] $r_i\ge0$ and $s_i\ge0$ for $1\le i\le m$.
\end{enumerate}
The sequence $((r_i,s_i))_{i=0}^m$ is the \emph{underlying walk} of $\lambda$.
If $m\ge1$ we put
\[
|\lambda|=\sum_{i=1}^m(\alpha_i+\beta_i),\qquad
\nu_\alpha(\lambda)=|\{i:\alpha_i>0\}|,\qquad
\nu_\beta(\lambda)=|\{i:\beta_i>0\}|.
\]
For $n\ge1$ let $\cL_n$ be the set of labeled  walks $\lambda$ with $m=0$ or
$|\lambda|\le n$, and for $\lambda\in\cL_n$ put
$d_n(\lambda)=\nu_\alpha(\lambda)+\nu_\beta(\lambda)-[\,m\ge1\text{ and }|\lambda|=n\,]$
(so $d_n(\lambda)=0$ if $m=0$).
\end{defn}

\begin{prop}\label{prop:walks}
For every $n\ge1$ there is a bijection $\Lambda_n\colon I_n(2143)\to\cL_n$ such
that $\des(\sigma)=d_n(\Lambda_n(\sigma))$ for all $\sigma\in I_n(2143)$.
Explicitly, let $(p_i,q_i,\rho_i)_{i=1}^k$ be the indexing of $\cE(\sigma)$ given
by Theorem~\ref{thm:fulton}(b), let $m=\lceil k/2\rceil$ and
$(p_0,q_0,\rho_0)=(0,n,0)$. Then $\Lambda_n(\sigma)=((u_i,\alpha_i,\beta_i))_{i=1}^m$
with
\begin{equation}\label{eq:labels}
u_i=(p_i-\rho_i)-(p_{i-1}-\rho_{i-1}),\qquad \alpha_i=p_i-p_{i-1},\qquad
\beta_i=q_{i-1}-q_i ,
\end{equation}
and the points of the underlying walk are
\begin{equation}\label{eq:points}
(r_i,s_i)=\bigl(\rho_i,\ n-p_i-q_i+\rho_i\bigr)\qquad(1\le i\le m),
\end{equation}
that is, the numbers of dots north-west and south-east of the essential box
$(p_i,q_i)$. Moreover, if $\sigma\neq\id_n$, then $m\ge1$,
\begin{equation}\label{eq:af}
f(\sigma)=\alpha_1\qquad\text{and}\qquad
a(\sigma)=\min\{\rho:\ (p,q,\rho)\in\cE(\sigma)\}=\min_{1\le i\le m}r_i ,
\end{equation}
so that $r(\sigma)=\alpha_1-\min_{1\le i\le m}r_i$. Moreover
\begin{equation}\label{eq:bl}
\ell(\sigma)=\beta_1,\qquad b(\sigma)=\min_{1\le i\le m}s_i,
\end{equation}
\begin{equation}\label{eq:pk}
\pk(\sigma)=\bigl|\{i\le m:\alpha_i\ge2\}\bigr|+\bigl|\{2\le i\le m:\beta_i\ge2\}\bigr|+[\,n-|\lambda|\ge2\,],
\end{equation}
where $\lambda=\Lambda_n(\sigma)$.
\end{prop}

\begin{proof}
\emph{Definition of $\Lambda_n$.} Let $\sigma\in I_n(2143)$ and index
$\cE(\sigma)=\{(p_i,q_i,\rho_i):1\le i\le k\}$ as in
Theorem~\ref{thm:fulton}(b). By Lemma~\ref{lem:inv-des}(a), $\cE(\sigma)$ is
stable under $(p,q,\rho)\mapsto(q,p,\rho)$. The sequence
$(q_{k+1-i},p_{k+1-i},\rho_{k+1-i})_{i=1}^k$ is an indexing of $\cE(\sigma)$
that satisfies \eqref{eq:F1}--\eqref{eq:F3}; by uniqueness of the
indexing,
\begin{equation}\label{eq:sym}
(p_i,q_i,\rho_i)=(q_{k+1-i},p_{k+1-i},\rho_{k+1-i})\qquad(1\le i\le k).
\end{equation}
Let $m=\lceil k/2\rceil$. If $i<k+1-i$, then by \eqref{eq:sym} and
\eqref{eq:F2}, $p_i-\rho_i<p_{k+1-i}-\rho_{k+1-i}=q_i-\rho_i$, so $p_i<q_i$;
if $i=k+1-i$ then $p_i=q_i$; and if $i>k+1-i$ then $p_i>q_i$. Hence the
triples of $\cE(\sigma)$ with $p\le q$ are those with $i\le m$, and
$p_m=q_m$ if and only if $k$ is odd. Define $\lambda=\Lambda_n(\sigma)$ by
\eqref{eq:labels}. Since
$\sum_{j\le i}\alpha_j=p_i$, $\sum_{j\le i}\beta_j=n-q_i$ and
$\sum_{j\le i}u_j=p_i-\rho_i$, we get \eqref{eq:points}. We check (S1) and (S2).
By \eqref{eq:F2}, $u_i\ge1$. By \eqref{eq:F1}, $\alpha_i,\beta_i\ge0$ for $i\ge2$,
while $\alpha_1=p_1>0$ and $\beta_1=n-q_1>0$ because $\Ess(\sigma)\subseteq[n-1]^2$.
Moreover
\[
\alpha_i+\beta_i-u_i=(q_{i-1}-\rho_{i-1})-(q_i-\rho_i),
\]
which is $\ge1$ by \eqref{eq:F2} if $i\ge2$, and equals $n-(q_1-\rho_1)\ge1$ if
$i=1$. Thus (S1) holds, and (S2) is \eqref{eq:F3}. Finally
$|\lambda|=p_m+n-q_m\le n$ because $p_m\le q_m$, so $\lambda\in\cL_n$.

\emph{The inverse map.} Let $\lambda=((u_i,\alpha_i,\beta_i))_{i=1}^m\in\cL_n$,
with underlying walk $((r_i,s_i))_{i=0}^m$. For $1\le i\le m$ put
\[
p_i=\sum_{j\le i}\alpha_j,\qquad q_i=n-\sum_{j\le i}\beta_j,\qquad \rho_i=r_i,
\qquad U_i=p_i-\rho_i=\sum_{j\le i}u_j,\qquad V_i=q_i-\rho_i .
\]
Put $k=2m-[\,m\ge1\text{ and }|\lambda|=n\,]$ and
\[
(p_i,q_i,\rho_i)=(q_{k+1-i},\ p_{k+1-i},\ \rho_{k+1-i})\qquad(m<i\le k).
\]
If $k=2m-1$, this is consistent for $i=m$ because $|\lambda|=p_m+n-q_m=n$ gives
$p_m=q_m$. By (S1), $1\le U_1<U_2<\dots<U_m$, and
$V_{i-1}-V_i=\alpha_i+\beta_i-u_i\ge1$ for $2\le i\le m$, so
$V_1>V_2>\dots>V_m$; moreover $V_m-U_m=q_m-p_m=n-|\lambda|\ge0$, with equality if
and only if $k=2m-1$. By construction \eqref{eq:sym} holds, and
\begin{align*}
(p_1,\dots,p_k)&=(p_1,\dots,p_m,\ q_{k-m},\dots,q_1),\\
(p_1-\rho_1,\dots,p_k-\rho_k)&=(U_1,\dots,U_m,\ V_{k-m},\dots,V_1).
\end{align*}
The first sequence is weakly increasing because $\alpha_i\ge0$, $\beta_i\ge0$ and
$p_m\le q_m\le q_{k-m}$. The second one is strictly increasing with
positive entries by the inequalities above, and because $U_m<V_m$ if
$k=2m$, while $U_m=V_m<V_{m-1}$ if $k=2m-1$. By \eqref{eq:sym}, the
sequences $(q_i)$ and $(q_i-\rho_i)$ are the reverses of these two sequences.
Finally, for $i\le m$ we have $\rho_i=r_i\ge0$ and
$p_i+q_i-\rho_i=n-s_i\le n$ by (S2), and the same holds for $i>m$
by \eqref{eq:sym}. Thus \eqref{eq:F1}--\eqref{eq:F3} hold, and
Theorem~\ref{thm:fulton}(b) provides a unique $\sigma\in S_n$ with
$\cE(\sigma)=\{(p_i,q_i,\rho_i)\}$, and $\sigma$ avoids $2143$. By
\eqref{eq:sym} and Lemma~\ref{lem:inv-des}(a),
$\cE(\sigma^{-1})=\cE(\sigma)$, so $\sigma^{-1}=\sigma$ by
Theorem~\ref{thm:fulton}(a). So $\sigma\in I_n(2143)$, and we set
$\Lambda_n'(\lambda)=\sigma$.

By construction, and since in both constructions the triples with $p\le q$
are exactly those with index $i\le m$, we have $\Lambda_n\circ\Lambda_n'=\mathrm{id}$;
and $\Lambda_n'\circ\Lambda_n=\mathrm{id}$ by Theorem~\ref{thm:fulton}(a),
because \eqref{eq:sym} recovers $\cE(\sigma)$ from its triples with $p\le q$.

\emph{Descents.} By Lemma~\ref{lem:inv-des}(b) and \eqref{eq:sym},
\[
\Des(\sigma)=\{p_1,\dots,p_k\}=\{p_1,\dots,p_m\}\cup\{q_1,\dots,q_m\},
\]
and $p_1\le\dots\le p_m\le q_m\le\dots\le q_1$. Since $p_0=0<p_1$ and
$q_1<n=q_0$, the number of distinct values among $p_1,\dots,p_m$ is
$\nu_\alpha(\lambda)$ and the number of distinct values among $q_1,\dots,q_m$ is
$\nu_\beta(\lambda)$. The two sets are disjoint if $p_m<q_m$, that is, if
$|\lambda|<n$, and they share exactly the element $p_m$ if $|\lambda|=n$. Hence
$\des(\sigma)=d_n(\Lambda_n(\sigma))$.

\emph{Initial fixed points and first descent.} Let $\sigma\neq\id_n$. Then
$\Des(\sigma)\neq\emptyset$, so $\Ess(\sigma)\neq\emptyset$ by
Lemma~\ref{lem:inv-des}(b), and $m\ge1$; by the previous paragraph the smallest
descent of $\sigma$ is $p_1=\alpha_1$. Let $a=a(\sigma)$; then $a<n$. If
$(p,q)\in D(\sigma)$ and $p\le a$, then $q<\sigma(p)=p\le a$, so
$\sigma^{-1}(q)=q<p$, contradicting $(p,q)\in D(\sigma)$; hence $p>a$, and $q>a$ by
Lemma~\ref{lem:inv-des}(a). Therefore the $a$ boxes $(h,h)$, $h\le a$, lie
north-west of every box of $D(\sigma)$, and $\rk_\sigma(p,q)\ge a$ for every
$(p,q)\in\Ess(\sigma)$. On the other hand, $\sigma(a+1)>a+1$, because
$\sigma(a+1)\neq a+1$ and the values $1,\dots,a$ are taken; since $\sigma$ is an
involution, also $\sigma^{-1}(a+1)>a+1$. Thus $(a+1,a+1)\in D(\sigma)$ and
$\rk_\sigma(a+1,a+1)=a$. By Cases~1 and~2 of the proof of
\cite[Lemma~3.10(a)]{Fulton1992}, every box of $D(\sigma)$ has an essential box
weakly south-east of it with the same rank. Hence the minimum of $\rk_\sigma$ on
$\Ess(\sigma)$ is $a$. By \eqref{eq:sym} this minimum is attained among the
triples with index $i\le m$, whose ranks are $\rho_i=r_i$ by \eqref{eq:points}.

\emph{Last run and final fixed points.} The largest descent is $q_1=n-\beta_1$,
so $\ell(\sigma)=\beta_1$. Let $b=b(\sigma)$ and $c=n-b$. If $(p,q)\in D(\sigma)$ and
$p>c$, then $\sigma(p)=p$, and $\sigma^{-1}(q)>p>c$ gives $\sigma^{-1}(q)=q>p$, while
$q<\sigma(p)=p$, a contradiction; so $p\le c$, and $q\le c$ by
Lemma~\ref{lem:inv-des}(a). Hence the $b$ final fixed points lie south-east of every
box of $D(\sigma)$, and $n-p-q+\rk_\sigma(p,q)\ge b$ for all $(p,q)\in D(\sigma)$.
Conversely, let $d=\sigma(c)$; then $d<c$ and $\sigma(d)=c$. Let $q^*$ be the largest
$j\le c-1$ with $\sigma(j)>d$ (it exists, since $\sigma(d)=c>d$). Then
$(d,q^*)\in D(\sigma)$, and the dots south-east of $(d,q^*)$ are exactly the $b$
final fixed points: a value $j\in(q^*,c-1]$ has $\sigma(j)\le d$, and the value $c$
lies in row $d$. Since every box of $D(\sigma)$ has an essential box weakly
south-east of it with the same rank \cite[Lemma~3.10(a)]{Fulton1992}, and
$n-p-q+\rho$ decreases weakly when $(p,q)$ moves south-east with $\rho$ fixed, the
minimum of $n-p-q+\rho$ over $\cE(\sigma)$ is $b$. This quantity is invariant under
\eqref{eq:sym}, so the minimum is attained at some $i\le m$, where it equals $s_i$
by \eqref{eq:points}.

\emph{Peaks.} A descent $i$ of $\sigma$ is a peak if and only if $i\ge2$ and
$i-1\notin\Des(\sigma)$, that is, if and only if the increasing run ending at $i$ has
length at least $2$. By the paragraph on descents, the runs ending at a descent
are: for each $i\le m$ with $\alpha_i>0$, the run ending at $p_i$, of length
$p_i-p_{i-1}=\alpha_i$; if $|\lambda|<n$, the run ending at $q_m$, of length
$q_m-p_m=n-|\lambda|$; for each $2\le i\le m$ with $\beta_i>0$, the run ending at
$q_{i-1}$, of length $q_{i-1}-q_i=\beta_i$ (if $|\lambda|=n$, then $p_m=q_m$ and this
also covers $i=m$). The last run, of length $\beta_1$, does not end at a descent.
This proves \eqref{eq:pk}.
\end{proof}

\begin{rem}\label{rem:endpoint}
The starting point $(r_0,s_0)=(0,0)$ of the underlying walk is the point that
\eqref{eq:points} associates with the virtual triple $(p_0,q_0,\rho_0)=(0,n,0)$: the
box $(0,n)$ has no dots north-west or south-east of it. The walk then follows the
triples of $\cE(\sigma)$ with $p\le q$, and it stops at the triple $(p_m,q_m,\rho_m)$,
whose point $(r_m,s_m)$ can be any point of $\N^2$. The triples with index $i>m$ are
not encoded, because they are determined by the first ones through \eqref{eq:sym}:
by \eqref{eq:points}, the point associated with $(p_i,q_i,\rho_i)$, $i>m$, is
$(r_{k+1-i},s_{k+1-i})$. Thus, for involutions, the walk starts at the origin but
its endpoint is free, and the only global condition on it is $|\lambda|\le n$.
If one encoded the whole essential set, adding the virtual triple
$(p_{k+1},q_{k+1},\rho_{k+1})=(n,0,0)$, whose box has no dots north-west or south-east
of it either, one would obtain an excursion from $(0,0)$ to $(0,0)$, namely the walk
of $\Lambda_n(\sigma)$ followed by its reversal. This folding is what makes the
involution case solvable: in Sections~\ref{sec:walks} and~\ref{sec:solve} only walks
with a prescribed starting point and a free endpoint occur. For general vexillary
permutations no folding is available, and the walk is an excursion
(Appendix~\ref{sec:vex}).
\end{rem}

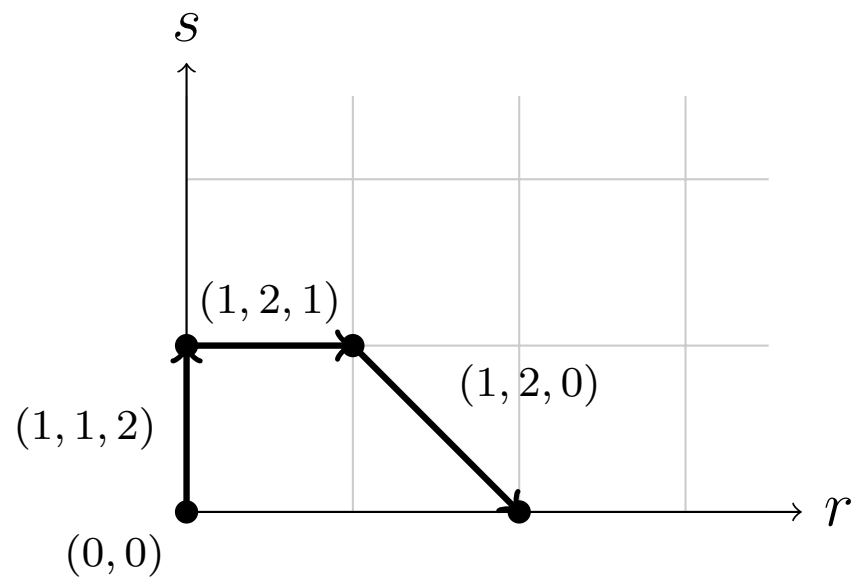
\begin{figure}[ht]
\centering
\begin{tikzpicture}[scale=1.1]
  \draw[step=1,gray!40] (0,0) grid (3.5,2.5);
  \draw[->] (0,0) -- (3.7,0) node[right] {$r$};
  \draw[->] (0,0) -- (0,2.7) node[above] {$s$};
  \draw[very thick,->] (0,0) -- (0,1);
  \draw[very thick,->] (0,1) -- (1,1);
  \draw[very thick,->] (1,1) -- (2,0);
  \foreach \p in {(0,0),(0,1),(1,1),(2,0)} \fill \p circle (0.07);
  \node[left]  at (-0.05,0.5) {\scriptsize $(1,1,2)$};
  \node[above] at (0.5,1) {\scriptsize $(1,2,1)$};
  \node[above right] at (1.5,0.5) {\scriptsize $(1,2,0)$};
  \node[below left] at (0,0) {\scriptsize $(0,0)$};
\end{tikzpicture}
\caption{The labeled  walk $\Lambda_8(\sigma)$ of $\sigma=7\,2\,6\,4\,8\,3\,1\,5$
(Example~\ref{ex:walk}). Its points $(0,1),(1,1),(2,0)$ are the numbers of dots
north-west and south-east of the essential boxes $(1,6),(3,5),(5,5)$; each step
is labeled by $(u_i,\alpha_i,\beta_i)$. The product
$\prod_i t^{[\alpha_i>0]+[\beta_i>0]}x^{\alpha_i+\beta_i}=t^2x^3\cdot t^2x^3\cdot t\,x^2=t^5x^8$
equals $t\cdot t^{\des(\sigma)}x^8$, since $|\lambda|=8=n$ (see the proof of
Proposition~\ref{prop:G2143W}).}
\label{fig:walk}
\end{figure}

\begin{ex}[Example~\ref{ex:essential}, continued]\label{ex:walk}
For $\sigma=7\,2\,6\,4\,8\,3\,1\,5$ the triples of $\cE(\sigma)$ with
$p\le q$ are $(1,6,0)$, $(3,5,1)$, $(5,5,2)$, so that $m=3$ and
\[
\Lambda_8(\sigma)=\bigl((u_i,\alpha_i,\beta_i)\bigr)_{i=1}^{3}
=\bigl((1,1,2),\ (1,2,1),\ (1,2,0)\bigr),
\]
with underlying walk $(0,0),(0,1),(1,1),(2,0)$ (Figure~\ref{fig:walk}). Here
$|\lambda|=8=n$, in accordance with the presence of the essential box
$(5,5)$ on the main diagonal, and
$d_8(\Lambda_8(\sigma))=3+2-1=4=\des(\sigma)$. Moreover $f(\sigma)=1=\alpha_1$ and
$a(\sigma)=0=\min(0,1,2)$.
\end{ex}

\begin{ex}
For $\sigma=1\,5\,7\,6\,2\,4\,3\in I_7(2143)$ we have
$\cE(\sigma)=\{(3,6,2),(4,4,1),(6,3,2)\}$ and
$\Lambda_7(\sigma)=((1,3,1),(2,1,2))$, with underlying walk $(0,0),(2,0),(1,0)$.
Indeed $\Des(\sigma)=\{3,4,6\}$, $f(\sigma)=3=\alpha_1$, $a(\sigma)=1=\min(r_1,r_2)$
and $r(\sigma)=2$. Note that the minimum is not attained at the first essential box.
\end{ex}

\section{Summing over the labels}\label{sec:walks}

Let $\lambda$ be a labeled walk and let $(\delta_i,\varepsilon_i)=(r_i-r_{i-1},s_i-s_{i-1})$
be the $i$-th step of its underlying walk. Since $\alpha_i=u_i+\delta_i$ and
$\beta_i=u_i+\varepsilon_i$, the label $(u_i,\alpha_i,\beta_i)$ is determined by the
step and by $u_i$. Define $\chi\colon\N\to R$ by $\chi(0)=1$, $\chi(1)=tx$ and
$\chi(k)=t\,\omega\,x^k$ for $k\ge2$. As we shall see in the proof of
Proposition~\ref{prop:G2143W}, the contribution of $\lambda$ to the generating
function of Theorem~\ref{thm:general} is essentially $\prod_i\chi(\alpha_i)\chi(\beta_i)$;
summing over all the labels compatible with a given step leads to the following
definition. For $(\delta,\varepsilon)\in\Z^2$ we define the \emph{weight} of
$(\delta,\varepsilon)$ as
\begin{multline}\label{eq:omega-def}
\wt(\delta,\varepsilon)=x^{|\delta|+|\varepsilon|}\Bigl(\omega^2t^2y+\omega t\,[\delta\varepsilon<0]
+\omega(1-\omega)t^2x^2\bigl([\delta\varepsilon<0]+[\delta=0]+[\varepsilon=0]\bigr)\Bigr)\\
+(1-\omega)^2t^2x^2\,[\delta=\varepsilon=0]\ \in R ;
\end{multline}
for $\omega=1$ it reduces to $x^{|\delta|+|\varepsilon|}(t^2y+t[\delta\varepsilon<0])$.
The next lemma justifies this definition. 

\begin{lem}\label{lem:omega}
For all $(\delta,\varepsilon)\in\Z^2$,
\[
\sum_{(u,\alpha,\beta)}\chi(\alpha)\,\chi(\beta)=\wt(\delta,\varepsilon),
\]
where the sum runs over the triples $(u,\alpha,\beta)\in\Z^3$ satisfying (S1)
with $\alpha-u=\delta$ and $\beta-u=\varepsilon$.
\end{lem}

\begin{proof}
Such a triple is determined by $u$, since $\alpha=u+\delta$ and
$\beta=u+\varepsilon$; condition (S1) becomes $u\ge1$, $u\ge-\delta$,
$u\ge-\varepsilon$ and $u\ge1-\delta-\varepsilon$, and $\alpha+\beta=2u+\delta+\varepsilon$.
Write $\chi(k)=\chi^\circ(k)+(1-\omega)tx\,[k=1]$, where $\chi^\circ(0)=1$ and
$\chi^\circ(k)=\omega t x^k$ for $k\ge1$. Recall that $y=\sum_{u\ge1}x^{2u}$, so that
$\sum_{u\ge u_0}x^{2u}=x^{2u_0-2}y$ for $u_0\ge1$. Since the sums below are symmetric
in $\delta$ and $\varepsilon$, it suffices to consider four cases to compute
$\sum\chi^\circ(\alpha)\chi^\circ(\beta)$.
\begin{itemize}
\item $\delta\ge0$, $\varepsilon\ge0$: the admissible values are $u\ge1$ and the sum
is $\omega^2t^2x^{\delta+\varepsilon}y$.
\item $\delta<0<\varepsilon$: the admissible values are $u\ge-\delta$; the term
$u=-\delta$ equals $\omega t\,x^{-\delta+\varepsilon}$, and the terms $u\ge1-\delta$
sum to $\omega^2t^2x^{-\delta+\varepsilon}y$.
\item $\delta<0$, $\varepsilon=0$: the admissible values are $u\ge1-\delta$ and the
sum is $\omega^2t^2x^{-\delta}y$.
\item $\delta<0$, $\varepsilon<0$: the admissible values are
$u\ge1-\delta-\varepsilon$ and the sum is $\omega^2t^2x^{-\delta-\varepsilon}y$.
\end{itemize}
Hence $\sum\chi^\circ(\alpha)\chi^\circ(\beta)=x^{|\delta|+|\varepsilon|}(\omega^2t^2y+\omega t[\delta\varepsilon<0])$.
Next, $\alpha=1$ means $u=1-\delta$, which is admissible if and only if
$\delta\le0\le\varepsilon$; then $\beta=1+|\delta|+|\varepsilon|\ge1$ and
$\chi^\circ(\beta)=\omega t\,x^{1+|\delta|+|\varepsilon|}$. Symmetrically, $\beta=1$ is
admissible if and only if $\varepsilon\le0\le\delta$, and $\alpha=\beta=1$ if and only if
$\delta=\varepsilon=0$ and $u=1$. Therefore
\begin{multline*}
\sum\chi(\alpha)\chi(\beta)=\sum\chi^\circ(\alpha)\chi^\circ(\beta)
+\omega(1-\omega)t^2x^{2+|\delta|+|\varepsilon|}\bigl([\delta\le0\le\varepsilon]+[\varepsilon\le0\le\delta]\bigr)\\
+(1-\omega)^2t^2x^2[\delta=\varepsilon=0],
\end{multline*}
and $[\delta\le0\le\varepsilon]+[\varepsilon\le0\le\delta]=[\delta\varepsilon<0]+[\delta=0]+[\varepsilon=0]$.
\end{proof}

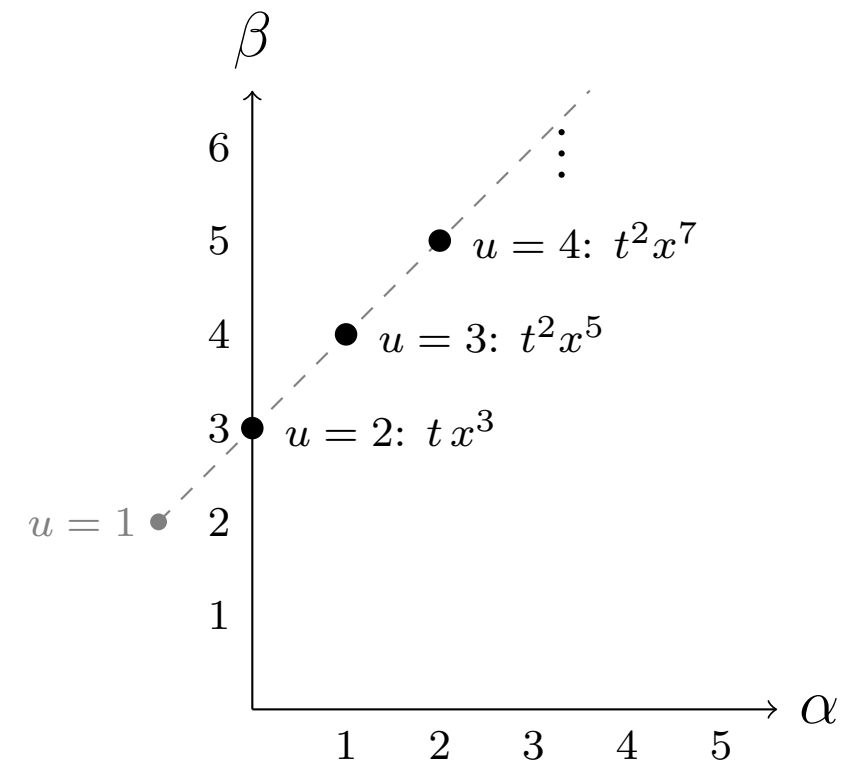
\begin{figure}[ht]
\centering
\begin{tikzpicture}[scale=0.62]
  \draw[->] (0,0) -- (5.6,0) node[right] {$\alpha$};
  \draw[->] (0,0) -- (0,6.6) node[above] {$\beta$};
  \foreach \k in {1,...,5} \node[below] at (\k,0) {\scriptsize $\k$};
  \foreach \k in {1,...,6} \node[left] at (0,\k) {\scriptsize $\k$};
  \draw[dashed,gray] (-1,2) -- (3.6,6.6);
  \fill[gray] (-1,2) circle (0.09);
  \node[left,gray] at (-1,2) {\scriptsize $u=1$};
  \fill (0,3) circle (0.12) node[right=2pt] {\scriptsize $u=2$: $t\,x^{3}$};
  \fill (1,4) circle (0.12) node[right=2pt] {\scriptsize $u=3$: $t^{2}x^{5}$};
  \fill (2,5) circle (0.12) node[right=2pt] {\scriptsize $u=4$: $t^{2}x^{7}$};
  \node at (3.3,6.1) {$\vdots$};
\end{tikzpicture}
\caption{The labels $(u,\alpha,\beta)=(u,u+\delta,u+\varepsilon)$ of the step
$(\delta,\varepsilon)=(-2,1)$: the admissible values are
$u\ge2$ (the value $u=1$ would give $\alpha<0$), each with weight
$t^{[\alpha>0]+[\beta>0]}x^{\alpha+\beta}$. Their sum is
$\wt(-2,1)=x^{3}(t+t^2y)$, in accordance with Lemma~\ref{lem:omega} (case $\omega=1$).}
\label{fig:step}
\end{figure}

See  Figure~\ref{fig:step} for an illustration of the previous lemma. 

A \emph{walk} of length $m\ge0$ starting at $(r,s)\in\N^2$ is a sequence
$\eta=((r_i,s_i))_{i=0}^m$ of points of $\N^2$ with $(r_0,s_0)=(r,s)$; there
is no restriction on the steps $(r_i-r_{i-1},s_i-s_{i-1})\in\Z^2$. Its
\emph{weight} is the product of the weights of its steps,
\[
\Wt(\eta)=\prod_{i=1}^m\wt(r_i-r_{i-1},\,s_i-s_{i-1}) .
\]
By \eqref{eq:omega-def}, $\Wt(\eta)\in t^mx^{L}R$, where
$L=\sum_i(|r_i-r_{i-1}|+|s_i-s_{i-1}|)$; since there are finitely many walks
with given starting point, $m$ and $L$, the series
\begin{equation}\label{eq:W}
W(r,s)=\sum_{m\ge0}\ \sum_{\eta}\Wt(\eta)\ \in R,
\end{equation}
where $\eta$ runs over the walks of length $m$ starting at $(r,s)$, is well
defined.

\begin{prop}\label{prop:G2143W}
For $d,e\in\Z$ put $\Delta W(d,e)=W(d,e)-W(d-1,e)-W(d,e-1)+W(d-1,e-1)$, where
$W(d,e)=0$ if $d<0$ or $e<0$. Then the left-hand side of Theorem~\ref{thm:general}
equals
\[
t\Bigl(1+tx+\frac{t\omega x^2}{1-x}\Bigr)\sum_{d,e\ge0}\Bigl(\omega\,\frac{w\bar wx^2}{1-w\bar wx^2}\cdot\frac{(xw)^d}{1-xv}
+(1-\omega)\,w\bar wx^2\,[d=0]\Bigr)\frac{(x\bar w)^e}{1-x\bar v}\,\Delta W(d,e).
\]
\end{prop}

\begin{proof}
By Proposition~\ref{prop:walks}, $\sigma\mapsto\Lambda_n(\sigma)$ maps
$I_n(2143)\setminus\{\id_n\}$ onto the nonempty labeled walks in $\cL_n$, and by
\eqref{eq:af}, \eqref{eq:bl} the statistics $a,r,b,\bar r$ depend only on the labeled
walk, not on $n$. A labeled walk $\lambda$ with $m\ge1$ belongs to $\cL_n$ exactly for
$n\ge|\lambda|$, with $\des(\sigma)=\nu_\alpha+\nu_\beta-[n=|\lambda|]$. By \eqref{eq:pk},
assigning to the $i$-th label the weight $\chi(\alpha_i)\chi(\beta_i)$ for $i\ge2$, and
$\chi(\alpha_1)\,t\,x^{\beta_1}$ for $i=1$ (the last run gives no peak), the
contributions of $n=|\lambda|$, $n=|\lambda|+1$ and $n\ge|\lambda|+2$ are the product of
these weights times $t^{-1}$, $x$ and $\omega x^{n-|\lambda|}$ respectively; their sum is
$\frac1t\bigl(1+tx+\frac{t\omega x^2}{1-x}\bigr)\,\chi(\alpha_1)\,t\,x^{\beta_1}\prod_{i\ge2}\chi(\alpha_i)\chi(\beta_i)$.

Write the underlying walk of $\lambda$ as a first step from $(0,0)$ to $(r_1,s_1)$
followed by a walk $\eta$ starting at $(r_1,s_1)$. The labeled walks with a given
underlying walk are obtained by choosing, independently for each $i$, a triple
$(u_i,\alpha_i,\beta_i)$ satisfying (S1) with $\alpha_i-u_i=\delta_i$,
$\beta_i-u_i=\varepsilon_i$; summing over the labels of the steps of $\eta$ gives
$\Wt(\eta)$, by Lemma~\ref{lem:omega}. In the first step $u_1\ge1$ is arbitrary,
$\alpha_1=u_1+r_1$ and $\beta_1=u_1+s_1$. By \eqref{eq:af} and \eqref{eq:bl},
$a(\sigma)=k$, $r(\sigma)=u_1+r_1-k$, $b(\sigma)=l$ and $\bar r(\sigma)=u_1+s_1-l$, where
$k$ and $l$ are the minima of the first and second coordinates of the points of
$\eta$. Translation by $(-k,-l)$ is a weight-preserving bijection from the walks
starting at $(r_1,s_1)$ whose points lie in $[k,\infty)\times[l,\infty)$ onto the walks
starting at $(r_1-k,s_1-l)$; by inclusion--exclusion the walks with minima exactly
$k$ and $l$ have total weight $\Delta W(r_1-k,s_1-l)$. Finally
$\chi(\alpha_1)=\omega t x^{\alpha_1}+(1-\omega)tx\,[\alpha_1=1]$, and $\alpha_1=1$ if and
only if $u_1=1$ and $r_1=k=0$. Summing over $u_1\ge1$, $k,l\ge0$ and
$d=r_1-k$, $e=s_1-l\ge0$ gives the statement.
\end{proof}

\section{Solution of the walk model}\label{sec:solve}

Let $\cF$ be the set of all functions $G\colon\N^2\to R$, and $\one\in\cF$
the constant function $1$. For $G\in\cF$ define $TG\in\cF$ by
\[
(TG)(r,s)=\sum_{(r',s')\in\N^2}\wt(r'-r,\,s'-s)\,G(r',s') .
\]
The sum converges in $R$ because
$\wt(r'-r,s'-s)\in t\,x^{|r'-r|+|s'-s|}R$ by \eqref{eq:omega-def}; for the
same reason $T(t^j\cF)\subseteq t^{j+1}\cF$ for all $j\ge0$.

\begin{lem}\label{lem:unique}
The function $W$ of \eqref{eq:W} is the unique $G\in\cF$ such that
$G=\one+TG$.
\end{lem}

\begin{proof}
Splitting a walk of length $m\ge1$ into its first step and the remaining walk
of length $m-1$ gives $W=\one+TW$. If $G=\one+TG$, then $G-W=T(G-W)$, hence
$G-W\in t^j\cF$ for every $j$ by induction, and $G=W$.
\end{proof}

Let $\theta_\omega$ be the series defined by \eqref{eq:starom}, and define
\begin{equation}\label{eq:gammaA}
g(r)=1-\frac{\theta_\omega-x}{1-x}\,\theta_\omega^{\,r}\qquad(r\in\N).
\end{equation}

\begin{prop}\label{prop:W}
For all $(r,s)\in\N^2$,
\[
W(r,s)=\frac{(1-x)^3\,g(r)\,g(s)}{\bigl(1-x+tx-t(1-\omega)x^2\bigr)\bigl((1-x)^2(1-tx)-t\omega x^2(3-x)\bigr)}.
\]
\end{prop}

\begin{proof}
By Lemma~\ref{lem:unique} it suffices to prove that the function $G(r,s)$ defined
by the right-hand side satisfies $G=\one+TG$. Since $G$ is a constant multiple of
$g\otimes g\colon(r,s)\mapsto g(r)g(s)$ and $T$ is linear, it suffices to compute
$T(g\otimes g)$. For $r\in\N$ put
\[
g^+(r)=\sum_{j>r}x^{j-r}g(j),\qquad
g^-(r)=\sum_{0\le j<r}x^{r-j}g(j),\qquad
g^{\pm}(r)=g^-(r)+g(r)+g^+(r).
\]
Since $[\delta\varepsilon<0]=[\delta>0][\varepsilon<0]+[\delta<0][\varepsilon>0]$,
\eqref{eq:omega-def} gives
\begin{multline*}
T(g\otimes g)(r,s)=\omega^2t^2y\,g^{\pm}(r)g^{\pm}(s)
+\bigl(\omega t+\omega(1-\omega)t^2x^2\bigr)\bigl(g^+(r)g^-(s)+g^-(r)g^+(s)\bigr)\\
+\omega(1-\omega)t^2x^2\bigl(g(r)g^{\pm}(s)+g^{\pm}(r)g(s)\bigr)+(1-\omega)^2t^2x^2g(r)g(s).
\end{multline*}
The sums $g^+(r)$, $g^-(r)$ are geometric; computing them for symbolic $r$ and
substituting the explicit expression of $\theta_\omega$, a computer algebra system
shows that $(TG)(r,s)=G(r,s)-1$ for all $(r,s)\in\N^2$.
\end{proof}

We are now in position to prove Theorems~\ref{thm:general} and~\ref{thm:refined} and Corollary~\ref{cor:peaks}.

\begin{proof}[Proof of Theorem~\ref{thm:general}]
Let $g(-1)=0$. For a formal variable $\xi$,
\[
\sum_{d\ge0}\xi^d\bigl(g(d)-g(d-1)\bigr)=(1-\xi)\sum_{d\ge0}\xi^dg(d)
=1-\frac{(\theta_\omega-x)(1-\xi)}{(1-x)(1-\xi\theta_\omega)}
=\frac{(1-\theta_\omega)(1-x\xi)}{(1-x)(1-\xi\theta_\omega)},
\]
and in particular $g(0)=(1-\theta_\omega)/(1-x)$. By Proposition~\ref{prop:W},
$\Delta W(d,e)$ is the constant of Proposition~\ref{prop:W} times
$\bigl(g(d)-g(d-1)\bigr)\bigl(g(e)-g(e-1)\bigr)$. Inserting this in
Proposition~\ref{prop:G2143W} and using the identity
\[
t\Bigl(1+tx+\frac{t\omega x^2}{1-x}\Bigr)
\frac{(1-x)^3(1-\theta_\omega)^2}{\bigl(1-x+tx-t(1-\omega)x^2\bigr)\bigl((1-x)^2(1-tx)-t\omega x^2(3-x)\bigr)}
=\frac{t\,\theta_\omega(1-x)^2}{x\,q_\omega},
\]
which is checked with a computer algebra system, gives the statement.
\end{proof}


\begin{proof}[Proof of Theorem~\ref{thm:refined}]
By Theorem~\ref{thm:general} with $\omega=\bar v=\bar w=1$, the generating function
of $I_n(2143)\setminus\{\id_n\}$ equals
$\frac{t\,x\,(1+x)\,\theta\,w}{(1-x\theta)(1-xv)(1-xw\theta)}$, which is
$\frac{w(\theta-x)}{(1-xv)(1-xw\theta)}$ by \eqref{eq:star}. The statement follows from
Proposition~\ref{prop:G3412}.
\end{proof}

\begin{proof}[Proof of Corollary~\ref{cor:peaks}]
The identity permutations have no descents and no peaks. Hence, by
Theorem~\ref{thm:general} with $v=w=\bar v=\bar w=1$,
\[
F_\omega=\frac1{1-x}+\frac{t\,x\,(1+x)\,\theta_\omega}{q_\omega(1-x\theta_\omega)}
\Bigl(1-\omega+\frac{\omega}{(1-x)(1-x\theta_\omega)}\Bigr),
\]
and a computer algebra system shows that $t\omega(F_\omega-1)=\theta_\omega/x-1$.
\end{proof}

\section{Monotone blocks and symmetric fillings}\label{sec:wilf-blocks}

In this section and in the next two, following \cite{Krattenthaler2006},
in every permutation plot positions increase from left to right and
values increase from bottom to top: the point $(p,\pi(p))$ is drawn
with Cartesian coordinates, and an increasing subsequence moves right
and up. This differs from the convention used for $D(\sigma)$ in
Section~\ref{sec:ess}, which follows \cite{Fulton1992}.
The auxiliary Ferrers boards below are indexed by position-block and
value-block numbers, both increasing in these directions. 

\begin{defn}\label{def:wilf-block}
For $m,n\in\N$, a blocking of $\pi\in S_n$ into $m$ blocks is a weak
composition $\ell_1+\cdots+\ell_m=n$; set
$b_0=0$, $b_i=\ell_1+\cdots+\ell_i$, and
$B_i=\{b_{i-1}+1,\ldots,b_i\}$. These are intervals of positions,
possibly empty. We apply the same intervals to the values in $[n]$.
The blocking is increasing, respectively decreasing, if the entries
of $\pi$ are increasing, respectively decreasing, within any
$B_i$. When $m=0$, only $n=0$ is allowed.
\end{defn}

\begin{lem}\label{lem:wilf-bars}
Let $G(x,t)=\sum_{n\ge0}x^n\sum_{\pi\in\mathcal A_n}t^{d(\pi)}$,
where $\mathcal A_n\subseteq S_n$ is a family of permutations sets with $\mathcal{A}_0$ containing the empty permutation. For increasing blockings put
$d=\des$, and for decreasing blockings put $d=\asc$.
If $H(x,u)$ enumerates the corresponding blocked permutations by
length and number of blocks, then
\begin{equation}\label{eq:wilf-bars}
 H(x,u)=\frac{1+u\bigl(G(x/(1-u),u)-1\bigr)}{1-u}.
\end{equation}
\end{lem}
\begin{proof}
Consider a given $\pi \in \mathcal A_n.$
For $n>0$, a blocking assigns weakly increasing labels in $[m]$ to
the $n$ positions. They rise strictly at each of the $d$ comparisons
incompatible with the block order. The number of choices is
$\binom{n+m-1-d}{n}$, with value zero when $m\le d.$ This adds a term equal to  $x^n u^{d+1}/(1-u)^{n+1}$ to the generating function $H.$ The empty permutation admits
one blocking for every $m\ge0,$ thus adds a contribute equal to $\frac{1}{1-u}$ to $H.$
\end{proof}

Fix one of the two block conditions in Definition~\ref{def:wilf-block}:
the entries of the permutation must increase within every position
block, or they must decrease within every position block. For a
permutation $\pi$ with such a blocking, define
\[
a_{ij}=\bigl|\{p\in B_i:\pi(p)\in B_j\}\bigr|,
\qquad A(\pi)=(a_{ij})_{1\le i,j\le m}.
\]
If $\pi$ is an involution, then $a_{ij}=a_{ji}$, since inversion
exchanges the points with position in $B_i$ and value in $B_j$
with those having position in $B_j$ and value in $B_i$.

Conversely, let $A=(a_{ij})$ be a symmetric $m\times m$ matrix
of nonnegative integers. Define
$\ell_i=\sum_{j=1}^{m}a_{ij}$ and the intervals $B_i$ as in
Definition~\ref{def:wilf-block}. Partition the \emph{positions}
in each $B_i$ into consecutive subintervals $P_{ij}$ of
cardinalities $a_{ij}$. If the blocks are required to be increasing,
arrange these subintervals from left to right in the order
$P_{i1},\ldots,P_{im}$. If the blocks are required to be
decreasing, arrange them from left to right in the order
$P_{im},\ldots,P_{i1}$. Thus $P_{ij}$ is the set of positions
in $B_i$ whose images must belong to the \emph{value} interval $B_j$.
For $p\in P_{ij}$, define $\pi(p)$ to be the element of
$P_{ji}$ with the same rank when the blocks are increasing,
and with the opposite rank when the blocks are decreasing.
The two subintervals have the same cardinality because $A$
is symmetric. In either case, this rule defines an involution:
the assignment from $P_{ji}$ back to $P_{ij}$ is its inverse.
Moreover, the entries in every position block have the
specified monotonicity, and $A(\pi)=A$. Conversely, that
monotonicity forces both the order of the subintervals
$P_{ij}$ and the assignment within each pair, so $\pi$ is unique.
More generally, if a blocked permutation $\sigma$ has matrix
$A(\sigma)$, then its inverse $\sigma^{-1}$, equipped with the
same blocking of positions and values, has matrix
$A(\sigma^{-1})=A(\sigma)^{\mathsf T}$. The row sums of $A$
are the values $\ell_1,\ldots,\ell_m$ of the block lengths,
and the length $n$ of $\pi$ is
$n=|A|=\sum_{i,j}a_{ij}$.

\begin{figure}[ht]
\centering
\begin{tikzpicture}[font=\small]
  \node at (0,1.0) {$A=(a_{ij})$};
  \node at (0,0) {$\displaystyle
    \begin{pmatrix}
      2&1&0\\
      1&0&1\\
      0&1&1
    \end{pmatrix}$};
  \node at (0,-1.15) {$\ell_1=3,\ \ell_2=2,\ \ell_3=2$};

  \draw[->,thick] (1.75,0)--(3.05,0);

  \begin{scope}[shift={(3.9,-1.86)},x=0.58cm,y=0.58cm]
    \draw[step=1,gray!35] (0,0) grid (7,7);
    \draw[blue,thick] (3,0)--(3,7) (5,0)--(5,7);
    \draw[blue,thick] (0,3)--(7,3) (0,5)--(7,5);
    \foreach \p/\v in {1/1,2/2,3/4,4/3,5/6,6/5,7/7}
      \fill (\p-0.5,\v-0.5) circle[radius=2.3pt];
    \foreach \i in {1,...,7} {
      \node[below] at (\i-0.5,0) {\i};
      \node[left] at (0,\i-0.5) {\i};
    }
    \node[below] at (3.5,-0.65) {positions};
    \node[rotate=90] at (-1.3,3.5) {values};
  \end{scope}
  \node at (5.93,-3.15) {$\pi=124\,|\,36\,|\,57=1243657$};
\end{tikzpicture}
\caption{The symmetric matrix $A$ and its associated increasing-blocked
involution. The row sums $(3,2,2)$ determine the common position and
value blocks $B_1=\{1,2,3\}$, $B_2=\{4,5\}$ and $B_3=\{6,7\}$,
whose boundaries are shown in blue. In particular, $a_{11}=2$
corresponds to the points $(1,1)$ and $(2,2)$. Positions increase
to the right and values increase upwards.}
\label{fig:wilf-matrix-involution}
\end{figure}
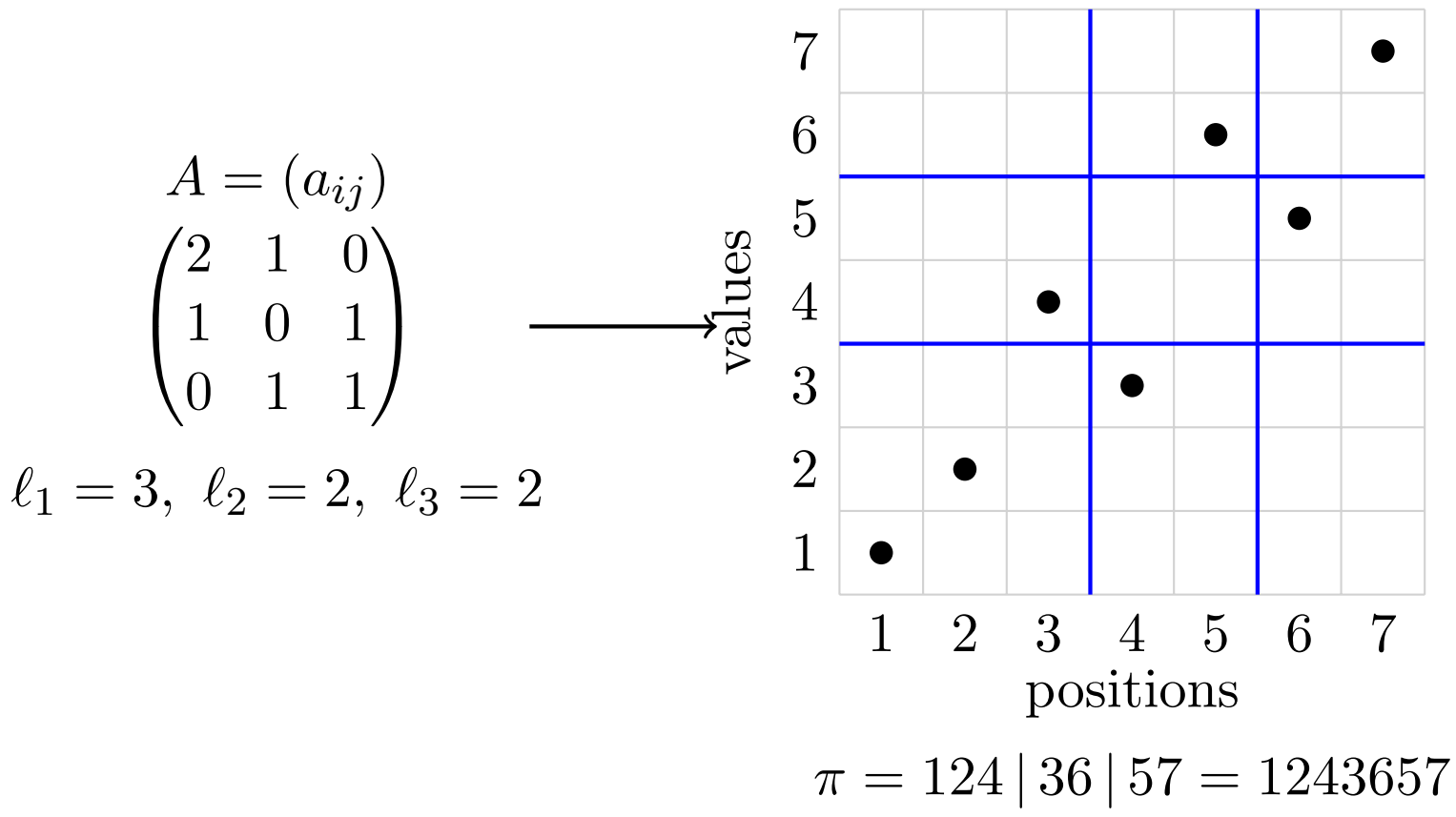

\begin{defn}\label{def:wilf-board}
For cells $(i,j),(p,q)\in[m]^2$, write $(p,q)\preceq(i,j)$ if
$p\le i$ and $q\le j$. A \emph{south-west order ideal} is
$\mathcal B\subseteq[m]^2$ such that
$(i,j)\in\mathcal B$ and $(p,q)\preceq(i,j)$ imply
$(p,q)\in\mathcal B$. A \emph{filling} $M$ assigns a
nonnegative integer to each cell of $\mathcal B$ and zero to every
other cell. An \emph{inside corner} is a maximal cell of
$\mathcal B$ under $\preceq$, and an \emph{outside corner} is a
minimal cell of $[m]^2\setminus\mathcal B$. A \emph{strict increasing
chain} consists of positive cells with both indices strictly
increasing; a \emph{strict decreasing chain} consists of positive
cells with first indices increasing and second indices decreasing.
A south-west rectangle is $[i]\times[j].$
\end{defn}
The directions chosen in the previous definition refer
to the plot convention above; in particular, an increasing
subsequence corresponds to a strict increasing chain in a
decreasing-blocked permutation.

\begin{figure}[ht]
\centering
\begin{tikzpicture}[x=0.78cm,y=0.78cm,font=\small]
  \fill[blue!9] (0,0) rectangle (5,2);
  \fill[blue!9] (0,2) rectangle (4,3);
  \fill[blue!9] (0,3) rectangle (2,4);
  \fill[blue!9] (0,4) rectangle (1,5);
  \draw[gray!55,step=1] (0,0) grid (5,5);
  \draw[blue!65!black,very thick]
    (5,0)--(5,2)--(4,2)--(4,3)--(2,3)--(2,4)--(1,4)--(1,5)--(0,5);

  \foreach \i/\j/\v in {1/1/1,2/1/2,4/1/1,1/2/1,
                        3/2/3,5/2/1,2/3/2,4/3/1,1/4/1,2/4/2,1/5/1}
    \node at (\i-0.5,\j-0.5) {$\v$};

  \foreach \i/\j in {1/1,3/2,4/3}
    \draw[red!80!black,very thick]
      (\i-0.5,\j-0.5) circle[radius=0.32];

  \foreach \i/\j in {5/2,2/3,1/4}
    \draw[orange!85!black,very thick]
      (\i-0.82,\j-0.82) rectangle (\i-0.18,\j-0.18);

  \foreach \i/\j in {5/2,4/3,2/4,1/5}
    \fill[blue!65!black] (\i-0.13,\j-0.13) circle[radius=0.095];
  \foreach \i/\j in {5/3,3/4,2/5}
    \fill[white] (\i-0.5,\j-0.5) circle[radius=0.11];
  \foreach \i/\j in {5/3,3/4,2/5}
    \draw[black,thick] (\i-0.5,\j-0.5) circle[radius=0.11];

  \foreach \i in {1,...,5} {
    \node[below] at (\i-0.5,0) {\i};
    \node[left] at (0,\i-0.5) {\i};
  }
  \node[below] at (2.5,-0.65) {$i$};
  \node[rotate=90] at (-0.75,2.5) {$j$};
\end{tikzpicture}
\caption{A south-west order ideal $\mathcal B\subseteq[5]^2$
(shaded cells) with a filling $M$; unlabelled cells have value zero.
Small blue dots mark the inside corners
$(5,2),(4,3),(2,4),(1,5)$, while open black circles mark the
outside corners $(5,3),(3,4),(2,5)$.
Red circles highlight the strict increasing chain
$(1,1),(3,2),(4,3)$; orange squares highlight the strict
decreasing chain $(5,2),(2,3),(1,4)$.
The south-west rectangle $[3]\times[2]$ consists of the cells in the
first three columns and in the two lowest rows. Both indices increase
in the directions shown.}
\label{fig:wilf-northwest-ideal}
\end{figure}
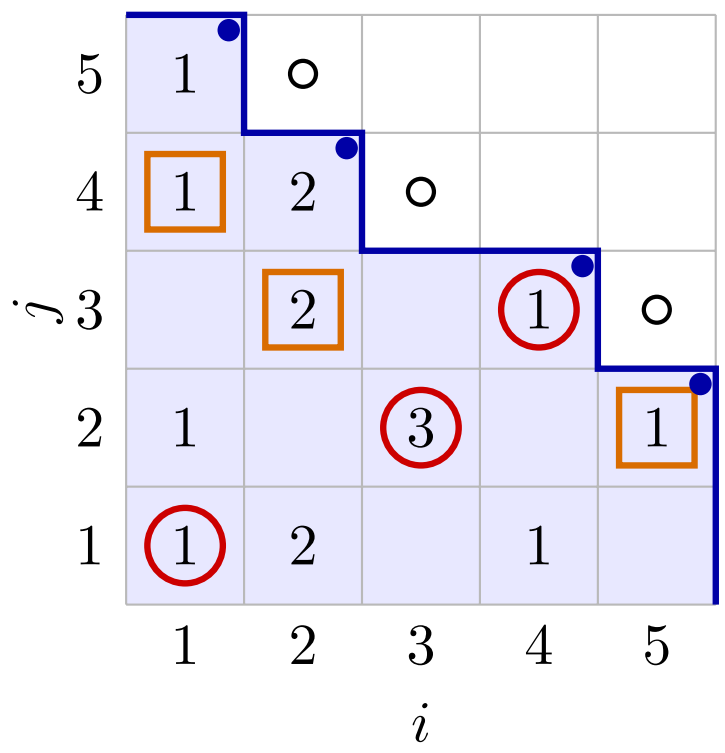

The inside corners of $\mathcal B$ determine it as the union of the
rectangles $[i]\times[j]$ ending at them. Its outside corners $C$
determine it by
\begin{equation}\label{eq:wilf-outside}
\mathcal B=[m]^2\setminus
\bigcup_{(i,j)\in C}\{i,\ldots,m\}\times\{j,\ldots,m\}.
\end{equation}
Either corner set is an antichain for $\preceq$. Both descriptions
include the empty board and the full square.

Let $\mathcal B\subseteq[m]^2$ be a south-west order ideal, the cell
$(i,j)$ being the unit square $[i-1,i]\times[j-1,j]$. Its \emph{boundary}
is the unique lattice path from $(0,m)$ to $(m,0)$, with unit steps to
the right, recorded by $+$, and down, recorded by $-$, such that the
cells of $\mathcal B$ are exactly the cells of $[m]^2$ lying below it;
it has $m$ steps of each kind. Let $v_0,\dots,v_{2m}$ be its vertices,
and for a vertex $v_k=(i,j)$ call $[i]\times[j]$ the south-west
rectangle at $v_k$. A \emph{boundary labelling} of $\mathcal B$ is a
sequence of partitions $\emptyset=\lambda^0,\lambda^1,\dots,\lambda^{2m}=\emptyset$
such that $\lambda^k/\lambda^{k-1}$ is a horizontal strip if the $k$-th
step is $+$, and $\lambda^{k-1}/\lambda^k$ is a horizontal strip if it
is $-$. If $\mathcal B$ is symmetric, a boundary labelling is
\emph{symmetric} if $\lambda^k=\lambda^{2m-k}$ for $0\le k\le 2m$.

\begin{lem}\label{lem:wilf-growth}
Let $\mathcal B\subseteq[m]^2$ be a symmetric south-west order ideal.
There are two bijections $\Gamma$ and $\Gamma'$ from the set of
symmetric fillings $M$ of $\mathcal B$ to the set of symmetric boundary
labellings of $\mathcal B$ such that, for $0\le k\le 2m$, the size of
$\lambda^k$ is the sum of the entries of $M$ in the south-west rectangle
at $v_k$, and the number of rows of $\lambda^k$ is the maximal length
of a strict decreasing chain (for $\Gamma$), respectively of a strict
increasing chain (for $\Gamma'$), of $M$ contained in that rectangle.
A symmetric boundary labelling is determined by
$\lambda^0,\dots,\lambda^m$, and the sum of the sizes of the strips
added or removed along the first $m$ steps is $|M|$.
\end{lem}
\begin{proof}
Let $\Gamma$ be the RSK growth diagram bijection of
\cite[Section~4.1]{Krattenthaler2006}, and let $\Gamma'$ be obtained from
the dual RSK$'$ growth diagram bijection of
\cite[Section~4.4]{Krattenthaler2006} by transposing all labels; dual
RSK$'$ labels differ by vertical strips, so the transposed labels differ
by horizontal strips. Our boards, boundaries and south-west rectangles
are the Ferrers shapes, D-R-sequences and rectangular regions to the left
and below of a corner of \cite{Krattenthaler2006}, and our strict
decreasing (respectively increasing) chains are its se-chains
(respectively ne-chains). Hence \cite[Theorems~7 and~8]{Krattenthaler2006},
respectively \cite[Theorems~11 and~12]{Krattenthaler2006}, give the
bijectivity for a fixed board and the assertions on sizes and numbers of
rows. By \cite[Section~5, discussion preceding Theorems~14 and~15]{Krattenthaler2006},
these bijections commute with the reflection in the main diagonal, hence
they restrict to bijections between symmetric fillings and symmetric
boundary labellings. Finally, along the complete boundary the strips
added at the steps $+$ and those removed at the steps $-$ both have total
size $|M|$; since the reflection exchanges the two halves of a symmetric
labelling and the two kinds of steps, the strips added or removed along
the first $m$ steps also have total size $|M|$.
\end{proof}

Write a partition with at most two rows as $(a,b)$, $a\ge b\ge0$.
For a function of these states define
\begin{equation}\label{eq:wilf-operators}
\begin{aligned}
(Uf)(a,b)&=\sum_{i\ge0}\sum_{j=0}^{a-b}x^{i+j}f(a+i,b+j),\\
(Df)(a,b)&=\sum_{i=0}^{a-b}\sum_{j=0}^{b}x^{i+j}f(a-i,b-j),\\
(D_1f)(a,b)&=\sum_{i=0}^{a-b}\sum_{j=0}^{b-1}x^{i+j}f(a-i,b-j).
\end{aligned}
\end{equation}
These ranges are exactly the possible additions and removals of
a horizontal strip; $D_1$ retains only target states with $b>0$.

Consider the following two systems of equations, for functions $H$, $L$,
$H'$ of the states $(a,b)$, and $J$ of the states $(a,b)$ with $b>0$,
with values in $\Q[[x,u]]$:
\begin{equation}\label{eq:wilf-3214-recur}
H=1+uUL+uDH,
\qquad L=x+uUL+ux^2DH,
\end{equation}
\begin{equation}\label{eq:wilf-1243-recur}
H'=1+uUH'+uD_1J,
\qquad
J=x+ux^2UH'+uD_1J\quad(b>0).
\end{equation}
Each of them has a unique solution, determined coefficientwise in $u$.

\begin{lem}\label{lem:wilf-common-solution}
Let $r\in x\Q[[x,u]]$ be the unique solution of
\begin{equation}\label{eq:wilf-kernel}
(r-x)(1-xr)=u(1+x)r,
\end{equation}
and let $(H,L)$ and $(H',J)$ be the solutions of
\eqref{eq:wilf-3214-recur} and \eqref{eq:wilf-1243-recur}. Then $H=H'$
at every state $(a,b)$, and $H(0,0)=r/x$.
\end{lem}
\begin{proof}
Put
\[
\gamma_1=\frac{(1-x)^2}{(1-x)^2-u(1+x)},\qquad
\gamma_2=-\frac{u(1+x)\,r\,\gamma_1}{(1-x)^2},\qquad
\gamma_3=-\frac{(r-x)\,\gamma_1}{1-x}.
\]
We claim that the solutions are
\begin{equation}\label{eq:wilf-solution}
\begin{aligned}
H&=H'=\gamma_1+\gamma_2r^{a-b}+\gamma_3(r^b-r^{a+1}),\\
L&=x\gamma_1+x\gamma_2r^{a-b}+xr\gamma_3(r^b-r^{a+1}),\\
J&=x\gamma_1+x\gamma_2r^{a-b}+(x/r)\gamma_3(r^b-r^{a+1})\qquad(b>0).
\end{aligned}
\end{equation}
These formulas were found by looking for solutions that are linear
combinations of the functions $(a,b)\mapsto p^aq^b$; since the solutions
are unique, it suffices to check that they satisfy
\eqref{eq:wilf-3214-recur} and \eqref{eq:wilf-1243-recur}. Applied to
these functions, the operators $U$, $D$ and $D_1$ give geometric sums;
computing them for symbolic $a,b$ and eliminating $u$ by means of
\eqref{eq:wilf-kernel}, a computer algebra system shows that both
systems are satisfied. Setting $a=b=0$ in the formula for $H$ gives
$H(0,0)=r/x$, again by \eqref{eq:wilf-kernel}.
\end{proof}

\section{Descents on $I_n(3214)$}\label{sec:wilf-3214}

\begin{prop}\label{prop:wilf-3214-board}
Increasing-blocked elements of $I_n(3214)$ with $m$ blocks are in
bijection with pairs $(\mathcal B,M)$ such that $\mathcal B$ is a
symmetric south-west order ideal in $[m]^2$, $M$ a symmetric filling
of $\mathcal B$, and no south-west rectangle in $\mathcal B$
contains a strict decreasing chain of length three.
If $C$ is the inside corner set, then $n=|M|+|C|$.
\end{prop}
\begin{proof}
Mark the right-to-left maxima of the permutation. They have decreasing
values, so each increasing position block contains at most one.
Inversion gives the same property for the value blocks. Their cells
form a symmetric antichain $C$. Every unmarked point is weakly
south-west of some marked point. Delete the marked points, giving a
symmetric filling $M$ of
$\mathcal B=\bigcup_{(i,j)\in C}[i]\times[j]$.

Conversely, add a unit at every inside corner and reconstruct the
increasing-blocked involution from its matrix. The added point is
last in its position block and largest in its value block, even
when that cell already contained points. The added points are exactly
the right-to-left maxima. An occurrence of $3214$ can be chosen
with its last point a right-to-left maximum, since one can replace
that point by a later larger one. The preceding three points are
unmarked and, by increasing order within position and value blocks,
form a strict decreasing chain of cells in its south-west rectangle.
The converse is immediate by adjoining the marked corner point.
\end{proof}

\begin{ex}
The involution $\pi=3412\in I_4(3214)$ with the increasing blocking
$34\,|\,12$ has matrix $A(\pi)$ with $a_{12}=a_{21}=2$ and
$a_{11}=a_{22}=0$. Its right-to-left maxima are $4$ and $2$, in the cells
$(1,2)$ and $(2,1)$. Hence $C=\{(1,2),(2,1)\}$,
$\mathcal B=\{(1,1),(1,2),(2,1)\}$, and $M$ has entries $1$ in $(1,2)$
and $(2,1)$; indeed $n=|M|+|C|=2+2$. For the involution $1243657$ of
Figure~\ref{fig:wilf-matrix-involution} the only right-to-left maximum
is $7$, so $C=\{(3,3)\}$, $\mathcal B=[3]^2$, and $M$ is obtained from
$A$ by decreasing $a_{33}$ by one.
\end{ex}

\begin{prop}\label{prop:wilf-3214-gf}
The descent generating function on $I(3214)$ is $F(x,t)$.
\end{prop}
\begin{proof}
Apply the bijection $\Gamma$ (ordinary RSK) of
Lemma~\ref{lem:wilf-growth}. Proposition~\ref{prop:wilf-3214-board}
restricts all boundary labels to two rows. Inside corners are
$+-$ turns. A turn in the interior of the first half contributes
$x^2$ for its reflected pair; a final $+$ gives a single central
corner of weight $x$. Each step has weight $u$ times $x$ to the
size of its strip. For a state $(a,b)$, let $H_-(a,b)$ be the generating
function of all the ways of completing the first half of a symmetric
boundary labelling from a vertex labelled $(a,b)$ that is the starting
vertex or is reached by a $-$ step, with the weights just described;
define $L_+(a,b)$ in the same way for a vertex reached by a $+$ step.
The blocked series is $H_-(0,0)$. Ending after
a $+$ gives the central factor $x$. Thus $(H_-,L_+)$ satisfies
\eqref{eq:wilf-3214-recur}, that is, it is the solution $(H,L)$ of that
system. Lemma~\ref{lem:wilf-common-solution} gives
$H_-(0,0)=\theta/[x(1-t)]$ after the auxiliary variables
$(x,u)$ are replaced by $(x(1-t),t)$.
Substituting in \eqref{eq:wilf-bars}, and using
$(1-x\theta)(\theta-x)=tx(1+x)\theta$ from \eqref{eq:star},
gives $(1+\theta)/(1-x\theta)=F(x,t)$.
\end{proof}

\section{Ascents on $I_n(1243)$}\label{sec:wilf-1243}

For a point $(i,\pi(i))$ let its \emph{increasing rank} be the
length of a longest increasing subsequence ending there. Call
the points of ranks one and two the \emph{core} $K$, and the
other points \emph{marked}. A marked point has two core points
strictly south-west of it forming an increasing pair.

\begin{lem}\label{lem:wilf-rank}
Fix a core $K$ arising from a permutation of $[n]$.
There is a unique completion avoiding $1243$ and a unique
completion avoiding $1234$. The correspondence between these
completions preserves involutions.
\end{lem}


\begin{proof}
The unused positions and values are called free. At a free position
$p$, a free value $v$ is \emph{eligible} if $(p,v)$ has a core
increasing pair strictly south-west of it. The marked points are the
members of the third basic subsequence with respect to $1234$
\cite[Definition~3.2.1]{West1990}; hence the completions described below
are the outputs of \cite[Algorithms~3.2.4 and~3.2.6]{West1990}, which fix
the core and are mutually inverse bijections between $S_n(1234)$ and
$S_n(1243)$ \cite[Theorem~3.2.3]{West1990}. A marked point has a core
increasing pair strictly south-west of it, so every completion assigns
eligible values to the free positions. In a $1234$-avoiding completion
the marked values decrease, since an increasing pair of marked points
extends a core pair to $1234$; so this completion is unique. In a
$1243$-avoiding completion, processing the free positions from left to
right, each one receives the least unused eligible value: choosing a
larger value while a smaller eligible value remains unused creates
$1243$ with the pair witnessing that smaller value; so this completion
is unique too. Matrix transposition preserves increasing ranks and each
of the two avoidance conditions; uniqueness shows that both completions
commute with inversion.
\end{proof}

\begin{ex}\label{ex:wilf-rank-completion}
Let $n=10$ and fix the core
\[
K=\{(1,10),(2,6),(3,7),(6,2),(7,3),(10,1)\}.
\]
The free positions and free values are both
$\{4,5,8,9\}$. At positions $4$ and $5$, the eligible free
values are $8$ and $9$: the core points $(2,6),(3,7)$
form a witnessing increasing pair. At positions $8$ and $9$,
the values $4$ and $5$ are also eligible, witnessed by
$(6,2),(7,3)$.

Choosing the least unused eligible value at each free position
gives the $1243$-avoiding completion
\[
\begin{array}{c|rrrrrrrrrr}
i
 &1&2&3&4&5&6&7&8&9&10\\ \hline
\pi_{1243}(i)
 &10&6&7&8&9&2&3&4&5&1
\end{array}
\]
whereas choosing the largest unused eligible value gives the
$1234$-avoiding completion
\[
\begin{array}{c|rrrrrrrrrr}
i
 &1&2&3&4&5&6&7&8&9&10\\ \hline
\pi_{1234}(i)
 &10&6&7&9&8&2&3&5&4&1.
\end{array}
\]
In both cases the points of ranks one and two are precisely $K$.
The first completion has cycles
\[
(1\,10)(2\,6)(3\,7)(4\,8)(5\,9),
\]
and the second has cycles
\[
(1\,10)(2\,6)(3\,7)(4\,9)(5\,8).
\]
Thus the two distinct completions are both involutions, as
asserted in Lemma~\ref{lem:wilf-rank}.
\end{ex}

\begin{lem}\label{lem:wilf-1243-block}
The bijection of Lemma~\ref{lem:wilf-rank} maps decreasing-blocked
involutions in $I_n(1243)$ to decreasing-blocked involutions in
$I_n(1234)$ with at most one marked point per block.
\end{lem}
\begin{proof}
In a $1243$-avoider two marked entries in adjacent positions must
form an ascent. Otherwise the core pair witnessing the second entry
lies before the first, and these four entries form $1243$.
An adjacent core entry followed by a marked one is always an
ascent: if it were a descent, the pair witnessing the second
would also witness the first. An adjacent marked entry followed
by a core entry is always a descent, since an ascent would
extend a triple to the core entry. Thus a decreasing block can
contain a marked point only in its first position. The completion
keeps all core entries and all marked positions fixed. Its
core--core and mixed comparisons inside each block have the same
signs, proving both directions.
\end{proof}

\begin{ex}\label{ex:wilf-decreasing-block-boundary}
Consider the decreasing-blocked involution
\[
\pi=1\,|\,5\,|\,3\,|\,42\,|\,6
   =153426\in I_6(1243).
\]
Its core is
\[
K=\{(1,1),(2,5),(3,3),(5,2)\},
\]
so the marked positions are $4$ and $6$. Applying the
bijection of Lemma~\ref{lem:wilf-rank} gives
\[
\pi'=1\,|\,5\,|\,3\,|\,62\,|\,4
    =153624\in I_6(1234).
\]
Both permutations are involutions and retain the same core and
block boundaries. Their fourth blocks, $42$ and $62$, are
decreasing and each contains one marked point at its first
position. The fifth block contains the other marked point.

\end{ex}

\begin{prop}\label{prop:wilf-1243-board}
Decreasing-blocked involutions in $I_n(1243)$ with $m$ blocks
correspond to pairs $(\mathcal B,M)$ consisting of a symmetric
south-west order ideal and a symmetric filling with
\begin{enumerate}
\item no strict increasing chain of length three in $M$;
\item a strict increasing pair of positive cells in
$[i-1]\times[j-1]$ for every outside corner $(i,j)$.
\end{enumerate}
If $C$ is the outside corner set, then $n=|M|+|C|$.
\end{prop}
\begin{proof}
Apply Lemma~\ref{lem:wilf-1243-block}. In the resulting $1234$-avoider
the marked values decrease, and Lemma~\ref{lem:wilf-1243-block} places the marked positions
in distinct blocks. Inversion gives distinct value blocks. Their cells form a symmetric antichain
$C$ with one point per cell. Delete them to obtain the core
matrix $M$. Each marked cell $(i,j)$ has a core increasing pair
in $[i-1]\times[j-1]$. No core cell can lie weakly north-east
of $(i,j)$, because that pair would make it a core triple.
The support of $M$ thus lies in the board \eqref{eq:wilf-outside},
whose outside corners are precisely $C$.

Conversely, add a unit in every outside corner. Condition 2
gives rank at least three to every added point. A positive cell
of $\mathcal B$ cannot follow an outside corner in an increasing
chain; distinct outside corners are incomparable. Condition 1
therefore makes the added points exactly the rank-three points
and the old points exactly the core. The resulting involution
avoids $1234$ and has at most one marked point per block.
Lemma~\ref{lem:wilf-1243-block} gives the inverse construction.
\end{proof}

For example, $14523$ avoids $1243$, and its entries in positions
$3$ and $5$ have rank three. They decrease in value, but they
are not adjacent. The blocking $1\,|\,4\,|\,52\,|\,3$ has outside
corners $(3,4)$ and $(4,3)$, with core cells $(1,1),(2,3),(3,2)$.

\begin{prop}\label{prop:wilf-1243-gf}
The ascent generating function on $I(1243)$ is $N(x,t)$.
\end{prop}
\begin{proof}
Use the bijection $\Gamma'$ (transposed dual RSK$'$) of
Lemma~\ref{lem:wilf-growth}. The absence of a core increasing
triple restricts labels to two rows. Outside corners are $-+$
turns; the second row of the partition at such a turn must be
positive by Proposition~\ref{prop:wilf-1243-board}(2).
An internal turn in the first half contributes $x^2$ for its
reflected pair, and a final $-$ contributes $x$ for the central
corner.

Define $H_+(a,b)$ as the generating function of the completions of the
first half of a symmetric boundary labelling from a vertex labelled
$(a,b)$ that is the starting vertex or is reached by a $+$ step, and,
for $b>0$, $J_-(a,b)$ in the same way for a vertex reached by a $-$
step, with the weights described above. A removal
ending at $b=0$ cannot lead to an admissible completion: further
removals keep $b=0$, and turning or stopping requires $b>0$.
By \eqref{eq:wilf-operators}, $(H_+,J_-)$ satisfies
\eqref{eq:wilf-1243-recur}, that is, it is the solution $(H',J)$ of that
system.
Lemma~\ref{lem:wilf-common-solution} gives the same blocked
series as \eqref{eq:wilf-3214-recur}. Lemma~\ref{lem:wilf-bars}
then gives $N(x,t)=(1+\theta)/(1-x\theta)$.
\end{proof}

We are now in position to prove Theorem~\ref{wilf_equidist_des}.

\begin{proof}[Proof of Theorem~\ref{wilf_equidist_des}]
For $n\ge1$, $\des(\id_n)=0$ and $\asc(\Omega_n)=0$; for $n=0,$ $S_0$ contains only
the empty permutation. Hence setting $v=w=1$ in Theorem~\ref{thm:refined} gives
\[
F=\frac1{1-x}+\frac{\theta-x}{(1-x)(1-x\theta)}=N ,
\]
that is, $\des$ on $I_n(2143)$ and $\asc$ on $I_n(3412)$ are equidistributed for
every $n\ge0$.

The right-hand side equals
$(1+\theta)/(1-x\theta).$ 

Substituting the auxiliary variables $x\mapsto x(1-t)$
and $u\mapsto t$ in \eqref{eq:wilf-kernel} gives
\eqref{eq:star}, hence $r=\theta$. The proof of
Propositions~\ref{prop:wilf-3214-gf} and
\ref{prop:wilf-1243-gf} now also proves that 
\[
O(x,t)=P(x,t)=\frac{1+\theta}{1-x\theta}.
\]

The expression by radicals is obtained by solving
\eqref{eq:star}, choosing the root with zero constant term, a fact that can be checked with a computer algebra system.

Let $\mathsf{C}(v)=\sum_{j\ge0}\Cat(j)v^j$, so that
$\mathsf{C}(v)=1+v\,\mathsf{C}(v)^2$.
It is easy to verify, through a computer algebra system or by hand, that 
$$t(N-1)=\frac{c}{1-c}\,\mathsf{C}\biggl(\frac{yc}{(1-c)^2}\biggr).$$
Hence
\[
t(N-1)=\sum_{j\ge0}\Cat(j)\,y^j\,\frac{c^{\,j+1}}{(1-c)^{2j+1}}
=\sum_{j\ge0}\Cat(j)\,y^j\sum_{k\ge j}\binom{k+j}{2j}c^{\,k+1},
\]
and, since $c=tz$, for $n\ge1$ the coefficient of $t^{k}x^n$ in $F$ is
$[x^n]\,z^{k+1}\sum_{j=0}^{k}\binom{k+j}{k-j}\Cat(j)y^j$.
\end{proof}

\section{Wilf classification by descents for involutions avoiding a pattern of length four}
\label{sec:wilf-classification}

For $\tau\in S$, write
\[
I_{n,\des}(\tau;t)=\sum_{\pi\in I_n(\tau)}t^{\des(\pi)}.
\]
The polynomials $I_{n,\asc}(\tau;t)$ are defined analogously. Notice that $I_{n,\asc}(\tau;t)=t^{n-1}I_{n,\des}(\tau;t^{-1})$ for $n\ge1$.
Two patterns $\tau$ and $\rho$ are \emph{$\des$-Wilf-equivalent} (\emph{$\asc$-Wilf-equivalent}, \emph{$\des$-$\asc$-Wilf-equivalent}, respectively) on involutions if
$I_{n,\des}(\tau;t)=I_{n,\des}(\rho;t)$ ($I_{n,\asc}(\tau;t)=I_{n,\asc}(\rho;t),$ $I_{n,\des}(\tau;t)=I_{n,\asc}(\rho;t),$ respectively) for every $n\ge0$.
The equalities proved in Theorem~\ref{wilf_equidist_des}
complete these classifications for patterns of length four (the same classification for patterns of length three is trivial).

In fact all the other correspondences are the trivial ones implied by  elementary reasons, namely the following facts. 
\begin{itemize}
    \item Inversion maps avoidance
of $\tau$ to avoidance of $\tau^{-1}$. Since an involution equals
its inverse, $I_n(\tau)=I_n(\tau^{-1})$ as sets. \item The reverse complement of a permutation $\pi$ is the permutation $\operatorname{rc}(\pi)$ such that
\[
\operatorname{rc}(\pi)(i)=n+1-\pi(n+1-i).
\] Reverse
complementation maps involutions to involutions, maps avoidance of $\tau$ to
avoidance of $\operatorname{rc}(\tau)$, and preserves the number
of descents. Thus
$I_{n,\des}(\tau;t)=I_{n,\des}(\operatorname{rc}(\tau);t)$.
\item For the two monotone patterns, the RSK correspondence (see, e.g.,~\cite{Knuth1970PERMUTATIONSMA}) sends an
involution to a single standard Young tableau $P$ and satisfies
$\operatorname{Des}(\pi)=\operatorname{Des}(P)$. Transposing
$P$ interchanges avoidance of $1234$ and $4321$ and
complements the descent set in $[n-1]$. Consequently, for
$n\ge1$,
\[
I_{n,\des}(4321;t)=I_{n,\asc}(1234;t).
\]
\end{itemize}

There are no other correspondences, as shown, for $n=8,$ in Table~\ref{tab:des-wilf-classes}. 

\begin{table}[ht]
\centering
\small
\setlength{\tabcolsep}{5pt}
\begin{tabular}{@{}l l r@{}}
\hline
Patterns $\tau$ & Coefficients $([t^k]I_{8,\des}(\tau;t))_{k=0}^{7}$
  & $I_{8,\des}(\tau;1)$\\
\hline
$1324$ & $(1,7,36,83,111,66,16,1)$ & $321$\\
$1234$ & $(0,0,3,54,144,105,16,1)$ & $323$\\
$1243,2134,3412$ & $(1,7,36,85,111,66,16,1)$ & $323$\\
$1432,2143,3214$ & $(1,16,66,111,85,36,7,1)$ & $323$\\
$4321$ & $(1,16,105,144,54,3,0,0)$ & $323$\\
$4231$ & $(1,16,66,111,87,38,7,1)$ & $327$\\
$2431,3241,4132,4213$ & $(1,16,86,132,106,44,10,1)$ & $396$\\
$1342,1423,2314,3124$ & $(1,10,44,110,138,86,16,1)$ & $406$\\
$2341,4123$ & $(1,12,54,122,149,86,16,1)$ & $441$\\
$2413,3142$ & $(1,16,64,140,142,76,16,1)$ & $456$\\
$3421,4312$ & $(1,16,86,157,131,56,12,1)$ & $460$\\
\hline
\end{tabular}
\caption{The eleven des-Wilf classes on involutions for patterns
of length four. For each row the middle column lists the coefficients
of the descent polynomial at $n=8$, from $t^0$ to $t^7$.
The last column records the total number of avoiding involutions.}
\label{tab:des-wilf-classes}
\end{table}

\appendix
\section{Vexillary permutations}\label{sec:vex}

The encoding of Section~\ref{sec:ess} does not depend on the symmetry of
involutions in an essential way. Let $\sigma\in S_n(2143)$, let
$(p_i,q_i,\rho_i)_{i=1}^k$ be the indexing of $\cE(\sigma)$ given by
Theorem~\ref{thm:fulton}(b), and add the two virtual triples
$(p_0,q_0,\rho_0)=(0,n,0)$ and $(p_{k+1},q_{k+1},\rho_{k+1})=(n,0,0)$. For
$1\le i\le k+1$ define $u_i,\alpha_i,\beta_i$ by \eqref{eq:labels}, and let
$\Lambda(\sigma)=((u_i,\alpha_i,\beta_i))_{i=1}^{k+1}$.

\begin{prop}\label{prop:vex}
For every $n\ge1$, $\Lambda$ is a bijection from $S_n(2143)$ onto the set of
sequences $((u_i,\alpha_i,\beta_i))_{i=1}^{\ell}$, $\ell\ge1$, satisfying (S1) and
(S2) of Definition~\ref{def:walks}, such that $(r_\ell,s_\ell)=(0,0)$ and
$\sum_iu_i=n$. The points of the underlying walk, which is an excursion from
$(0,0)$ to $(0,0)$, are given by \eqref{eq:points}, and
\[
\des(\sigma)=\bigl|\{i<\ell:\ \alpha_i>0\}\bigr|,\qquad
\des(\sigma^{-1})=\bigl|\{i<\ell:\ \beta_i>0\}\bigr| .
\]
\end{prop}

\begin{proof}
The argument is that of Proposition~\ref{prop:walks}, without the symmetry
\eqref{eq:sym}. Conditions \eqref{eq:F1}--\eqref{eq:F3} translate into (S1)
and (S2) for the steps $1\le i\le k$ exactly as there. For the last step,
$u_{k+1}=n-(p_k-\rho_k)\ge1$ because $p_k\le n-1$,
$\alpha_{k+1}+\beta_{k+1}-u_{k+1}=q_k-\rho_k\ge1$ by \eqref{eq:F2}, and
$\alpha_{k+1}=n-p_k\ge0$, $\beta_{k+1}=q_k\ge0$; moreover $\sum_iu_i=n$ and the walk
ends at $(\rho_{k+1},n-p_{k+1}-q_{k+1}+\rho_{k+1})=(0,0)$. Conversely, given such a
sequence, the walk returns to $(0,0)$, so
$\sum_i\alpha_i=\sum_i\beta_i=\sum_iu_i=n$; setting
$p_i=\sum_{j\le i}\alpha_j$, $q_i=n-\sum_{j\le i}\beta_j$ and $\rho_i=r_i$ for
$i\le k=\ell-1$, conditions \eqref{eq:F1}--\eqref{eq:F3} follow from (S1) and (S2)
as in Proposition~\ref{prop:walks}; in particular
$q_k-\rho_k=n-\sum_{j\le k}(\alpha_j+\beta_j-u_j)=\alpha_\ell+\beta_\ell-u_\ell\ge1$.
Theorem~\ref{thm:fulton} gives the bijection. Finally, $\Des(\sigma)$ is the set of
rows of $\Ess(\sigma)$ and $\Des(\sigma^{-1})$ the set of its columns
(Lemma~\ref{lem:inv-des}), and the number of distinct values among
$p_1,\dots,p_k$ (resp.\ $q_1,\dots,q_k$) is the number of $i\le k$ with
$\alpha_i>0$ (resp.\ $\beta_i>0$).
\end{proof}

Summing over the labels as in Section~\ref{sec:walks}, now with weight
$x^{u}t^{[\alpha>0]}t'^{[\beta>0]}$ for each step except the last one, which
carries the weight $x^u$, one obtains
\[
\sum_{n\ge1}\ \sum_{\sigma\in S_n(2143)}t^{\des\sigma}\,t'^{\des\sigma^{-1}}x^n
=\frac{x}{1-x}\sum_{(r,s)\in\N^2}x^{r+s}\,Q(r,s),
\]
where $Q(r,s)$ is the total weight of the walks in $\N^2$ from $(0,0)$ to
$(r,s)$ (of any length $m\ge0$), each step $(\delta,\varepsilon)$ having weight
\[
\widetilde{\wt}(\delta,\varepsilon)
=x^{\delta^-+\varepsilon^-}\Bigl(\frac{tt'x}{1-x}+t'\,[\delta<0<\varepsilon]+t\,[\varepsilon<0<\delta]\Bigr),
\]
with $\delta^-=\max(0,-\delta)$; the factor $x^{r+s}\cdot x/(1-x)$ is the weight of
the last step, from $(r,s)$ to $(0,0)$. As in Lemma~\ref{lem:omega}, and with the
same proof, $\widetilde{\wt}(\delta,\varepsilon)$ is the sum of
$t^{[\alpha>0]}\,t'^{[\beta>0]}\,x^{u}$ over the triples $(u,\alpha,\beta)$
satisfying (S1) with $\alpha-u=\delta$ and $\beta-u=\varepsilon$.

We do not know how to solve this model. Here we explain why the method of
Section~\ref{sec:solve} does not apply. First, for involutions the symmetry \eqref{eq:sym} folds the essential set in two
(Remark~\ref{rem:endpoint}): only its first half is encoded, so that the walk starts
at the origin but its endpoint is free. Hence only the series $W(r,s)$ of all walks
starting at a given point is needed, and the equation $W=\one+TW$ has the product
solution of Proposition~\ref{prop:W}. For a general vexillary permutation there is no
such symmetry, and the walk runs over the whole essential set, from the virtual triple
$(0,n,0)$ to the virtual triple $(n,0,0)$. Since neither virtual box has dots
north-west or south-east of it, both correspond to the point $(0,0)$: the walk starts
and ends at the origin (Proposition~\ref{prop:vex}). One therefore needs the two-point
series $Q(r,s)$, and not
only the series of walks with a free endpoint. Second, the answer is not algebraic. Indeed, for $t=t'=1$ the
series above counts vexillary permutations, which are equinumerous with
$1234$-avoiding permutations \cite{BackelinWestXin2007} and whose generating function is known to be D-finite but not algebraic.
Hence no ansatz leading to an algebraic series, such as the product form
$Ag(r)g(s),$ where $A\in R,$  of Proposition~\ref{prop:W}, can work.
We leave the determination of
this generating function, which also refines the distribution of
$(\des\sigma,\des\sigma^{-1})$ on vexillary permutations, as an open problem. However, we are able to furnish a functional equation satisfied by $Q(u,v)=\sum_{r,s\ge0}Q(r,s)\,u^rv^s$. Note that $Q(x,v)$,
$Q(u,x)$ and $Q(x,x)$ are well defined, and that, by the formula above, the
generating function of vexillary permutations is $\frac{x}{1-x}\,Q(x,x)$.

\begin{prop}\label{prop:vexFE}
Let $T'=tt'x/(1-x)$ and
\[
K(u,v)=(1-u)(1-v)(u-x)(v-x)-tt'x(1-x)\,uv-t'xv(1-u)(v-x)-txu(1-v)(u-x).
\]
Then
\begin{multline*}
K(u,v)\,Q(u,v)=(1-u)(1-v)(u-x)(v-x)-t'xv(1-u)(v-x+tx)\,Q(x,v)\\
-txu(1-v)(u-x+t'x)\,Q(u,x)+T'x^2(1-u)(1-v)\,Q(x,x).
\end{multline*}
\end{prop}

\begin{proof}
For a series $G(u)=\sum_rg_ru^r$, summing $u^{r+\delta}$ over the steps $\delta\ge0$
gives $G(u)/(1-u)$, over the steps $\delta>0$ gives $uG(u)/(1-u)$, and summing
$x^{-\delta}u^{r+\delta}$ over the steps $\delta<0$ with $r+\delta\ge0$ gives
\[
\sum_r g_r\sum_{k=0}^{r-1}x^{r-k}u^k=\frac{x\,\bigl(G(u)-G(x)\bigr)}{u-x}=:N_uG .
\]
Splitting a walk ending at $(r,s)$ into its last step and the rest, and using
the formula for $\widetilde{\wt}$, we get
\[
Q=1+T'\,L_uL_vQ+t'\,\frac{v}{1-v}\,N_uQ+t\,\frac{u}{1-u}\,N_vQ,
\qquad L_u=\frac1{1-u}+N_u,
\]
where $N_v$ and $L_v$ are defined in the same way with respect to $v$.
Multiplying by $(1-u)(1-v)(u-x)(v-x)$ and collecting the terms in $Q(u,v)$,
$Q(x,v)$, $Q(u,x)$ and $Q(x,x)$ gives the stated equation; this last step is a
routine computation, checked with a computer algebra system. We have also
checked the equation coefficientwise, for the coefficients of $u^iv^jx^k$ with
$i,j\le3$ and $k\le5$.
\end{proof}

\section{Verification code}\label{app:code}
In this appendix we add the Mathematica code used to perform all the algebraic checks present in the paper.

\begin{lstlisting}        
ClearAll["Global`*"];
z = x/(1 - x); y = x^2/(1 - x^2); c = t z;
P = 1 + x^2 - t x (1 + x); D0 = Expand[P^2 - 4 x^2];
th = (P - Sqrt[D0])/(2 x);
enne = (1 - x^2 + x^2 t - x t -
    Sqrt[Expand[x^4 (t - 1)^2 + x^3 (2 t^2 - 2 t) + x^2 (t^2 - 2 t - 2) -
      2 t x + 1]])/(2 t x^2);
(* theta_omega of !\eqref{eq:starom}!; om stands for omega *)
qom = 1 - t (1 - om) x;
Bom = 1 + x^2 - t x - t om x^2 - t (1 - om) x^3;
D1 = Expand[Bom^2 - 4 x^2 qom^2];
thom = (Bom - Sqrt[D1])/(2 x qom);
(* formal zero test *)
zero[e_] := Module[{n},
   n = Numerator[Together[e /. Sqrt[a_] :> Sqrt[Expand[a]]]];
   n = n /. Power[a_, k_Rational] /; Denominator[k] == 2 :>
      rt[Expand[a]]^(2 k);
   n = Expand[n] //. rt[a_]^m_Integer /; m >= 2 :> a rt[a]^(m - 2);
   Expand[n] === 0];
(* Section !\ref{prel}!: theta and theta_omega *)
check1 = zero[x th^2 - P th + x];
check2 = SeriesCoefficient[th, {x, 0, 0}] === 0;
check3 = zero[x qom thom^2 - Bom thom + x qom] &&
   SeriesCoefficient[thom, {x, 0, 0}] === 0;
check4 = zero[(thom /. om -> 1) - th];
(* Proposition !\ref{prop:N}! *)
motz[n_] := n - 1 - x (1 + t (n - 1)) - x^2 n (1 + t (n - 1));
check5 = Together[-t x motz[1 + (T/x - 1)/t] - (x T^2 - P T + x)] === 0;
check6 = zero[motz[enne]] && SeriesCoefficient[enne, {x, 0, 0}] === 1;
check7 = zero[th - x (1 + t (enne - 1))];
check8 = zero[y (t (enne - 1))^2 - (1 - c) t (enne - 1) + c];
(* Proposition !\ref{prop:W}!: G = 1 + T G, with G(r,s) = Aom gg(r) gg(s),
   gg(r) = 1 - ga T^r, T = theta_omega, ga = (T - x)/(1 - x), r, s symbolic;
   T G is computed from !\eqref{eq:omega-def}! by formal geometric sums *)
Aom = (1 - x)^3/((1 - x + t x - t (1 - om) x^2)*
     ((1 - x)^2 (1 - t x) - t om x^2 (3 - x)));
gg[j_] := 1 - ga T^j;
gplus[r_] = Sum[x^(j - r) gg[j], {j, r + 1, Infinity},
   GenerateConditions -> False];
gminus[r_] = Sum[x^(r - j) gg[j], {j, 0, r - 1},
   GenerateConditions -> False];
gall[r_] = gminus[r] + gg[r] + gplus[r];
TG = om^2 t^2 y gall[r] gall[s] +
   (om t + om (1 - om) t^2 x^2) (gplus[r] gminus[s] + gminus[r] gplus[s]) +
   om (1 - om) t^2 x^2 (gg[r] gall[s] + gall[r] gg[s]) +
   (1 - om)^2 t^2 x^2 gg[r] gg[s];
res = Numerator[Together[
    Expand[PowerExpand[a0 TG + 1 - a0 gg[r] gg[s]]] //.
     {T^(r + s + k_.) :> a b T^k, T^(r + k_.) :> a T^k,
      T^(s + k_.) :> b T^k, x^(r + s + k_.) :> p q x^k,
      x^(r + k_.) :> p x^k, x^(s + k_.) :> q x^k}]];
check9 = FreeQ[res, r | s] &&
   And @@ (zero /@ (Last /@ CoefficientRules[res, {a, b, p, q}] /.
         {a0 -> Aom, ga -> (T - x)/(1 - x)} /. T -> thom));
(* Proof of Theorem !\ref{thm:general}! *)
check10 = Together[(1 - xi) Sum[xi^d gg[d], {d, 0, Infinity},
         GenerateConditions -> False] -
      (1 - T) (1 - x xi)/((1 - x) (1 - xi T)) /.
     ga -> (T - x)/(1 - x)] === 0;
mu = t (1 + t x + t om x^2/(1 - x));
check11 = zero[mu Aom (1 - thom)^2 - t thom (1 - x)^2/(x qom)];
Sd[xi_] := (1 - thom) (1 - x xi)/((1 - x) (1 - xi thom));
fromProp = mu Aom (om (w wb x^2/(1 - w wb x^2)) Sd[x w]/(1 - x v) +
     (1 - om) w wb x^2 (1 - thom)/(1 - x)) Sd[x wb]/(1 - x vb);
RHS = t x thom w wb (1 - wb x^2)/(qom (1 - x vb) (1 - x wb thom)) *
   (1 - om + om (1 - w x^2)/((1 - x v) (1 - x w thom) (1 - w wb x^2)));
check12 = zero[fromProp - RHS];
(* Proof of Theorem !\ref{thm:refined}! *)
check13 = zero[(RHS /. {om -> 1, vb -> 1, wb -> 1}) -
    w (th - x)/((1 - x v) (1 - x w th))];
(* Proof of Theorem !\ref{wilf_equidist_des}! *)
check14 = Together[1/(1 - x) + (T - x)/((1 - x) (1 - x T)) -
     (1 + T)/(1 - x T)] === 0;
check15 = zero[(1 + th)/(1 - x th) - enne];
(* Proof of Theorem !\ref{wilf_equidist_des}!: M = t(N-1)(1-c)/c satisfies
   M = 1 + vv M^2 with vv = y c/(1-c)^2, and M(0) = 1 *)
M = t (enne - 1) (1 - c)/c; vv = y c/(1 - c)^2;
check16 = zero[M - 1 - vv M^2] && SeriesCoefficient[M, {x, 0, 0}] === 1;
(* Proof of Corollary !\ref{cor:peaks}! *)
Fom = 1/(1 - x) + (RHS /. {v -> 1, w -> 1, vb -> 1, wb -> 1});
check17 = zero[t om (Fom - 1) - (thom/x - 1)];
(* Proposition !\ref{prop:vexFE}!: QA = Q(u,v), QB = Q(x,v), QC = Q(u,x),
   QD = Q(x,x) *)
Tp = t tp x/(1 - x);
NuQ = x (QA - QB)/(u - x); NvQ = x (QA - QC)/(v - x);
LuLvQ = QA/((1 - u) (1 - v)) + x (QA - QC)/((1 - u) (v - x)) +
   x (QA - QB)/((1 - v) (u - x)) +
   x^2 (QA - QB - QC + QD)/((u - x) (v - x));
opId = QA - (1 + Tp LuLvQ + tp v/(1 - v) NuQ + t u/(1 - u) NvQ);
Kuv = (1 - u) (1 - v) (u - x) (v - x) - t tp x (1 - x) u v -
   tp x v (1 - u) (v - x) - t x u (1 - v) (u - x);
check18 = Together[(1 - u) (1 - v) (u - x) (v - x) opId -
     (Kuv QA - ((1 - u) (1 - v) (u - x) (v - x) -
        tp x v (1 - u) (v - x + t x) QB -
        t x u (1 - v) (u - x + tp x) QC +
        Tp x^2 (1 - u) (1 - v) QD))] === 0;
(* Proof of Theorem !\ref{wilf_equidist_des}!: !\eqref{eq:wilf-kernel}! with
   x -> x(1-t), u -> t is !\eqref{eq:star}! *)
check19 = zero[(th - x (1 - t)) (1 - x (1 - t) th) - t (1 + x (1 - t)) th];
(* Lemma !\ref{lem:wilf-common-solution}!: rr stands for r, and u is
   eliminated by means of !\eqref{eq:wilf-kernel}!. A function
   (a,b) -> c p^a q^b is encoded as {c, p, q}; PA[p] = p^a and PB[q] = q^b
   are treated as independent indeterminates. opU, opD[.,1], opD[.,0] are
   U, D, D_1 of !\eqref{eq:wilf-operators}!, summed as geometric series *)
uu = (rr - x) (1 - x rr)/((1 + x) rr);
g1 = (1 - x)^2/((1 - x)^2 - uu (1 + x));
g2 = -rr uu (1 + x) g1/(1 - x)^2; g3 = -g1 (rr - x)/(1 - x);
ClearAll[PA, PB];
PA[1] = 1; PB[1] = 1;
PA[s_ w_] := PA[s] PA[w]; PB[s_ w_] := PB[s] PB[w];
PA[Power[s_, -1]] := 1/PA[s]; PB[Power[s_, -1]] := 1/PB[s];
ev[F_] := Total[(#[[1]] PA[#[[2]]] PB[#[[3]]]) & /@ F];
opU[F_] := Total[(#[[1]] PA[#[[2]]] PB[#[[3]]] *
      (1 - x #[[3]] PA[x #[[3]]]/PB[x #[[3]]])/
      ((1 - x #[[2]]) (1 - x #[[3]]))) & /@ F];
opD[F_, e_] := Total[(#[[1]] PA[#[[2]]] PB[#[[3]]] *
      (1 - (x/#[[2]]) PA[x/#[[2]]]/PB[x/#[[2]]]) *
      (1 - (x/#[[3]])^e PB[x/#[[3]]])/
      ((1 - x/#[[2]]) (1 - x/#[[3]]))) & /@ F];
Hs = {{g1, 1, 1}, {g2, rr, 1/rr}, {g3, 1, rr}, {-g3 rr, rr, 1}};
Ls = {{x g1, 1, 1}, {x g2, rr, 1/rr}, {x rr g3, 1, rr},
   {-x rr g3 rr, rr, 1}};
Js = {{x g1, 1, 1}, {x g2, rr, 1/rr}, {(x/rr) g3, 1, rr},
   {-(x/rr) g3 rr, rr, 1}};
z0[e_] := Expand[Numerator[Together[e]]] === 0;
check20 = And @@ (z0 /@ {ev[Hs] - 1 - uu opU[Ls] - uu opD[Hs, 1],
     ev[Ls] - x - uu opU[Ls] - uu x^2 opD[Hs, 1],
     ev[Hs] - 1 - uu opU[Hs] - uu opD[Js, 0],
     ev[Js] - x - uu x^2 opU[Hs] - uu opD[Js, 0]});
check21 = z0[(ev[Hs] /. {PA[_] -> 1, PB[_] -> 1}) - rr/x];
{check1, check2, check3, check4, check5, check6, check7, check8, check9,
 check10, check11, check12, check13, check14, check15, check16, check17,
 check18, check19, check20, check21}
\end{lstlisting}

\section*{Declarations}
We used Claude (Anthropic) to assist with some aspects
of preparing this manuscript (editing, reformulating propositions, and checking proofs). The author assumes responsibility for all content.

\addcontentsline{toc}{section}{Bibliography}
\bibliographystyle{abbrv}
\bibliography{references}

\end{document}